\documentclass[oneside, reqno] {amsart}

\usepackage{xypic}
\input xy
\xyoption{all}
\usepackage{epsfig}
\usepackage{amsthm}
\usepackage{amssymb}
\usepackage{amsmath}
\usepackage{amscd}
\usepackage{color}
\usepackage[T1]{fontenc}
\usepackage{upgreek}
\usepackage[font=scriptsize]{caption}
\usepackage{wrapfig}

\usepackage{graphicx}
\usepackage{subfigure}

\newcommand{\bg}{\begin{equation}}
\newcommand{\ed}{\end{equation}}
\newcommand{\bga}{\begin{eqnarray}}
\newcommand{\eda}{\end{eqnarray}}
\newcommand{\pf}{\textbf{Proof:\ }}

\def\cbdu{\par{\raggedleft$\Box$\par}}

\newtheorem {Theorem}  {Theorem}

\numberwithin{Theorem}{section}

\newtheorem {Lemma}[Theorem]  {Lemma}
\newtheorem {Proposition}[Theorem]{Proposition}
\theoremstyle{definition}

\theoremstyle{remark}

\expandafter\chardef\csname pre amssym.def
at\endcsname=\the\catcode`\@ \catcode`\@=11
\def\undefine#1{\let#1\undefined}
\def\newsymbol#1#2#3#4#5{\let\next@\relax
 \ifnum#2=\@ne\let\next@\msafam@\else
 \ifnum#2=\tw@\let\next@\msbfam@\fi\fi
 \mathchardef#1="#3\next@#4#5}
\def\mathhexbox@#1#2#3{\relax
 \ifmmode\mathpalette{}{\m@th\mathchar"#1#2#3}%
 \else\leavevmode\hbox{$\m@th\mathchar"#1#2#3$}\fi}
\def\hexnumber@#1{\ifcase#1 0\or 1\or 2\or 3\or 4\or 5\or 6\or 7\or 8\or
 9\or A\or B\or C\or D\or E\or F\fi}

\font\teneufm=eufm10 \font\seveneufm=eufm7 \font\fiveeufm=eufm5
\newfam\eufmfam
\textfont\eufmfam=\teneufm \scriptfont\eufmfam=\seveneufm
\scriptscriptfont\eufmfam=\fiveeufm

\catcode`\@=\csname pre amssym.def at\endcsname

\newcounter{remark}
\newcommand{\supp}{{\mathit supp}\,}

\newcommand{\e}{\epsilon}

\def  \12  {{\frac{1}{2}}}

\def\build#1_#2^#3{\mathrel{\mathop{\kern 0pt#1}\limits_{#2}^{#3}}}

\numberwithin{equation}{section}

\begin{document}

\title[Non-uniqueness of active scalar equations]{Non-uniqueness of forced active scalar equations with even drift operators}


\author [Mimi Dai]{Mimi Dai}

\address{Department of Mathematics, Statistics and Computer Science, University of Illinois at Chicago, Chicago, IL 60607, USA}
\email{mdai@uic.edu}

\author [Susan Friedlander]{Susan Friedlander}

\address{Department of Mathematics, University of Southern California, Los Angeles, CA 90089, USA}
\email{susanfri@usc.edu}

\thanks{M.D. is partially supported by the NSF grants DMS--2009422 and DMS--2308208. }

\begin{abstract}

We consider forced active scalar equations with even and homogeneous degree 0 drift operator on $\mathbb T^d$. Inspired by the non-uniqueness construction for dyadic fluid models \cite{DF,FK}, by implementing a sum-difference convex integration scheme we obtain non-unique weak solutions for the active scalar equation in space $C_t^0C_x^\alpha$ with $\alpha<\frac{1}{2d+1}$. We note that in 1D, the regularity $\alpha<\frac13$ is sharp as the energy identity is satisfied for solutions in $C^\alpha$ with $\alpha>\frac13$. Without external forcing, Isett and Vicol \cite{IV} constructed non-unique weak solutions for such active scalar equations with spatial regularity $C_x^\alpha$ for $\alpha<\frac{1}{4d+1}$.

\bigskip

KEY WORDS: forced active scalar equation; non-uniqueness; convex integration sum-difference scheme.

\hspace{0.02cm}CLASSIFICATION CODE: 35Q35, 35Q86, 76D03.
\end{abstract}

\maketitle

\section{Introduction}
\label{sec-int}

\subsection{Background}
The active scalar equation with external forcing
\begin{equation}\label{ase}
\begin{split}
\partial_t\theta+ u \cdot\nabla \theta=&-\nu\Lambda^\gamma \theta+f, \\
u=&\ T[\theta],\\
\nabla\cdot u=&\ 0
\end{split}
\end{equation}
describes a number of physical phenomena arising in fluid dynamics. The unknown $\theta$ is a real-valued scalar function, while $u$ is the drift velocity  defined from $\theta$ through the nonlocal Zygmund operator $T$. The given function $f$ denotes the external buoyancy forcing.  The parameter $\nu\geq 0$ is the dissipation coefficient, and $\gamma>0$ indicates the strength of the dissipation. The operator $T$ has Fourier symbol $m(\xi)$ which is even, homogeneous of degree 0, and satisfies $\xi\cdot m(\xi)=0$. We consider (\ref{ase}) on $\mathbb T^d\times[0,\infty)$. 

Particular physical examples of (\ref{ase}) with even drift operators include the incompressible porous media (IPM) equation \cite{Bear, CFO} and the magnetogeostrophic (MG) equation \cite{FV1, Moff, ML}. These physical models have attracted attention due to their application in various physics contexts and their connection to hydrodynamic equations. 

The class of active scalar equations with odd drift operators, including the surface quasi-geostrophic equation (SQG) \cite{CMT}, has also been extensively studied in the literature. The different symmetry features of the even and odd classes result in different ill/well-posedness theories. The cancellation property for the odd class of scalar equations is beneficial in establishing well-posedness, see \cite{CaV, CV, KNV, Mar, Res}; while such cancellation structure is absent for the even class. The main objective of this paper is to investigate the ill-posedness phenomena for (\ref{ase}) with even operators $T$ through the lens of convex integration techniques. 

A pair $(\theta, u)$ is a weak solution of (\ref{ase}) if the equations in (\ref{ase}) are satisfied in the distributional sense. In the inviscid case $\nu=0$, the existence of global weak solution for the active scalar equation with even operators is a challenging problem in the framework of classical energy method. 
Nevertheless, it was shown in \cite{Shv} that there are infinitely many bounded weak solutions for (\ref{ase}) with $\nu=0$ and $f\equiv 0$ via the application of convex integration techniques, which were developed for Euler equations in \cite{DLS1, DLS2}. In separate works \cite{CFG, Sze}, non-unique bounded weak solutions were also constructed for the 2D IPM equation based on the convex integration method.  In the recent work \cite{CCF} for the 2D IPM, the authors obtained infinitely many mixing solutions in Sobolev space by combining convex integration, contour dynamics and pseudodifferential operators techniques. Similar result for the 2D IPM with curved interfaces was established in \cite{CFM}.
We emphasize that the convex integration method in the aforementioned works is rooted in the Tartar framework through the concept of subsolution.   
In particular, the solutions constructed in \cite{CFG, Shv, Sze} are in the space $L^\infty_{t,x}$. 

In the time line of the progress toward solving Onsager's conjecture (verifying $\frac13$ H\"older regularity threshold for energy conservation) for the Euler equation, it was first shown the existence of bounded weak solutions violating the energy conservation in \cite{DLS2}, then improvements were obtained in \cite{CDLS, DLS1, DLS3, Is0, BDLIS, BDLS} by producing continuous and $C^\alpha$ dissipative solutions for $\alpha<\frac15$. Eventually dissipative solutions with spatial regularity $C^\alpha$ for any $\alpha<\frac13$ were constructed in \cite{Is}. In this development, the improvement from constructing bounded weak solutions to H\"older continuous solutions relies on a crucial cancellation property which involves the use of stationary plane wave solutions for the Euler equation. The benefit of taking such plane waves as building blocks is that interference terms between different waves can be controlled.

Coming back to active scalar equations with even drift operators, an analogous Onsager's conjecture is that $\frac13$ spatial regularity is expected to guarantee energy conservation for the solution. However, it is not trivial to adapt the ideas for the Euler equations in the hope of obtaining wild weak solutions that are continuous or $C^\alpha$ for $\alpha>0$. The obstacle is that interference terms in different waves for active scalar equations can not be controlled efficiently due to the lack of a similar cancellation structure as that for the Euler equation. Remarkably the authors of \cite{IV} discovered a new mechanism for producing cancellations between overlapping waves for active scalar equations with even operators, and constructed non-unique dissipative solutions in $C_x^{\alpha}$ with $\alpha<\frac{1}{4d+1}$ for (\ref{ase}) with $\nu=0$ and $f\equiv 0$. This mechanism can also be used to construct dissipative solutions for the Euler equation in $C_x^{\alpha}$ with $\alpha<\frac15$, which is consistent with the result for (\ref{ase}) in one dimension, that is, when $d=1$. The cancellation comes from the vanishing of self-interaction terms which is due to the property $\xi\cdot m(\xi)=0$ (the divergence free condition). 
Such a cancellation determines that the iteration convex integration scheme is essentially based on one dimensional oscillations. Hence considering the problem in $d$-dimension, i.e. on $\mathbb T^d$, requires $d$ stages to correct the stress error in the iteration step from $R_q$ to $R_{q+1}$. This explains why the spatial regularity obtained in \cite{IV} depends on dimension, which appears to be counter intuitive in some sense since one expects to have more flexibility in the higher dimensional case.

We note that the forced surface quasi-geostrophic (SQG) equation (an active scalar equation with odd drift operator) was recently studied in \cite{BHP2, DP, DP2}. For the forced stationary SQG, it was shown that there are more than one solutions in the space $C_x^{\alpha}$ with $\alpha<-\frac14$ in \cite{DP}. As a contrast, a non-trivial weak solution was constructed in $C_x^{\alpha}$ with $\alpha<-\frac13$ for stationary SQG without external forcing in \cite{CKL}. For the evolutionary SQG without external forcing, non-unique weak solutions with spatial regularity $C_x^\alpha$ for $\alpha<-\frac15$ were constructed independently in \cite{BSV} and \cite{IM}. While for the forced evolutionary SQG, the authors of \cite{BHP2} and \cite{DP2} recently constructed non-unique weak solutions with spatial regularity $C_x^\alpha$ for $\alpha<0$. The construction of \cite{BHP2} is in the framework of \cite{BSV} and the construction of \cite{DP2} is in the framework of \cite{IM}. Moreover in both \cite{DP} and \cite{DP2}, the authors exploit the flexibility due to the presence of forcing through the sum-difference formulation of two solutions for the underlying system, which was inspired by \cite{DF, FK}. Such sum-difference formulation will be adapted in the current paper as well. More details will be provided in Section \ref{sec-outline}.


Comparing the active scalar equations with even and odd drift operators, we observe that it seems much harder to construct weak solutions with higher regularity that violate uniqueness and the energy law in the odd case. The reason is that the cancellation property in the odd case presents an obstruction, see \cite{BSV}. 

Among the active scalar equations with even operator, the MG equation is of particular interest since the operator $T$ has an unbounded region in the Fourier space (c.f. \cite{FV2}). Due to the unboundedness and evenness of $T$, ill-posedness for (\ref{ase}) with $\nu=0$, $f\equiv 0$ and the MG operator $T$ was shown in \cite{FV2} in the sense that there is no Lipschitz solutions map at the initial time. While when $\nu>0$, the unforced MG equation was shown well-posed in \cite{FV1} since the diffusion term plays a dominant role. In the case of fractional diffusion for the MG equation, i.e. (\ref{ase}) with $\nu>0$ and $f\equiv 0$, the authors of \cite{FRV} identified a threshold value $\gamma=\frac12$ by proving that: the MG equation with $\gamma\in(\frac12, 1)$ is locally well-posed, the MG equation with $\gamma\in (0,\frac12)$ is ill-posed and the MG equation with $\gamma=\frac12$ is globally well-posed for small initial data. 


\subsection{Main result}
In this current paper we study the active scalar equation (\ref{ase}) with $\nu\geq 0$ and $f\not\equiv 0$.  The purpose is to test whether the flexibility of allowing an external forcing can lead to the construction of wild solutions that reach the critical spatial regularity $\frac13$ for the energy law. We adapt the cancellation mechanism discovered in \cite{IV} in our construction, and thus our result also depends on dimension. In particular, the wild solutions we obtain have spatial regularity $C_x^{\alpha}$ with $\alpha<\frac{1}{2d+1}$ and violate uniqueness. It is clear that this regularity is higher than $\frac{1}{4d+1}$ for the non-forced case in \cite{IV}. Moreover, in one dimension $d=1$, we have $\frac{1}{2d+1}=\frac13$ and hence the regularity of the non-unique solutions reaches the threshold of $\frac13$.



\medskip

The main result is stated below.
\begin{Theorem}\label{thm}
Let $\alpha<\frac{1}{2d+1}$, $0\leq \gamma<1-\alpha$ and $\zeta<\frac{1}{2d}$. There exists $f\in C_t^0C_x^{2\alpha-1}$ such that there are more than one solutions $\theta\in C_t^0C_x^{\alpha}(\mathbb T^d)\cap C_t^\zeta C_x^{0}(\mathbb T^d)$ to (\ref{ase}) with external forcing $f$.
\end{Theorem}

Theorem \ref{thm} implies that the forced MG equation with $\gamma<1-\frac{1}{2d+1}$ is ill-posed due to the lack of uniqueness.

The paper is organized as follows. Section \ref{sec-no-f} gives a heuristic analysis for the result of non-forced active scalar equation with non-odd drift operators which was proven in \cite{IV}; Section \ref{sec-outline} provides a heuristic analysis for the forced active scalar equation with an even operator. In Section \ref{sec-pre} we lay out technical preparations and the main iterative process. Section \ref{sec-iteration} is devoted to the proof of the main iteration statement. Section \ref{sec-proof} concludes the proof of Theorem \ref{thm}.



\bigskip

\section{Heuristics of non-uniqueness for (\ref{ase}) with $\nu=0$ and $f\equiv0$}
\label{sec-no-f}

We provide an outline of heuristics for the earlier result of Isett and Vicol in their article \cite{IV} concerning unforced active scalar equations.
We present the analysis in the latter notation of \cite{BV} for the Navier-Stokes equation.
\begin{Theorem}\label{thm-no-f}
Let $\alpha<\frac{1}{4d+1}$ and $I$ be an open time interval. There exist non-trivial solutions $\theta\in C^{\alpha}_{x,t}(\mathbb T^d\times \mathbb R)$ to (\ref{ase}) with $\nu=0$ and $f\equiv 0$ such that $\theta(x,t)=0$ for $t\notin I$.
\end{Theorem}

First, consider the approximating systems in 2D
\begin{equation}\label{eq-q}
\begin{split}
\partial_t\theta_q+\nabla\cdot(u_q\theta_q)=&\ \nabla\cdot\widetilde R_q,\\
u_q=&\ T[\theta_q].
\end{split}
\end{equation}
Assume the image of the even part of the multiplier $m$ contains $d$ linearly independent vectors given by
\begin{equation}\notag
\begin{split}
A_j=m(\xi_j)+m(-\xi_j), \ \ \ j=1, 2, ..., d, \ \ |\xi_j|=1.
\end{split}
\end{equation}
The stress vector $\widetilde R_q$ can be decomposed as 
\[\widetilde R_q=c_{1,q}A_1+c_{2,q}A_2=: c_{1,q}A_1+R_q.\]
Without loss of generality, we assume $|c_{1,q}|\geq |c_{2,q}|$. The goal is to construct a new solution such that the principal part $c_{1,q}A_1$ in the stress error gets reduced.

We specify the index $I=(k, \pm)\in \mathbb Z\times \{\pm\}:= \Omega$. Denote $\bar I=(k, \mp)$. For $I\in\Omega$, 
let $\theta_{I,q+1}$ and $\xi_{I}$ be the amplitude and phase functions respectively, satisfying
\[\theta_{\bar I}=\bar \theta_I, \ \ \xi_{\bar I}=-\xi_I.\]  
Moreover, $\xi_I$ is advected by $u_q$ on a short time interval $\tau_q$ with initial state $\widehat \xi_I$.
The increment $\Theta_{q+1}=\theta_{q+1}-\theta_q$ is constructed to take the form
\begin{equation}\notag
\begin{split}
\Theta_{I,q+1}=&\ \mathbb P_{I,\lambda_{q+1}}(e^{i\lambda_{q+1}\xi_{I}}\theta_{I, q+1})=e^{i\lambda_{q+1}\xi_{I}}(\theta_{I, q+1}+\delta\theta_{I,q+1}),\\
\Theta_{q+1}=&\sum_{I\in\Omega} \Theta_{I,q+1}
\end{split}
\end{equation}
where the error term $\delta\theta_{I,q+1}$ comes from the application of a microlocal lemma.
Applying the microloccal lemma again yields
\begin{equation}\notag
\begin{split}
U_{I,q+1}=&\ T[\Theta_{I,q+1}]= e^{i\lambda_{q+1}\xi_I}(u_{I,q+1}+\delta u_{I}) \ \ \mbox{with} \ \
u_{I,q+1}= m(\nabla\xi_{I})\theta_{I,q+1}\\
U_{q+1}=&\ T[\Theta_{q+1}]=\sum_{I\in\Omega}T[\Theta_{I,q+1}].
\end{split}
\end{equation}
We also have 
\[u_{q+1}=u_q+U_{q+1}=T[\theta_q]+ T[\Theta_{q+1}].\]
The tuple $(\theta_{q+1}, u_{q+1}, \widetilde R_{q+1})$ is a solution of (\ref{eq-q}) with $q$ replaced by $q+1$ and 
with the new stress error $\widetilde R_{q+1}$ satisfying
\begin{equation}\label{R-q1}
\begin{split}
\nabla\cdot \widetilde R_{q+1}=& \left(\partial_t+u_{q}\cdot\nabla \right)\Theta_{q+1}+\nabla\cdot(U_{q+1}\theta_q)\\
&+ \nabla\cdot\sum_{J\neq \bar I}U_{J,q+1}\Theta_{I, q+1}\\
&+ \nabla\cdot\sum_{I\in\Omega}\left(U_{I,q+1}\Theta_{\bar I, q+1}+c_{1,q}A_1+R_q\right)\\
=&: \nabla\cdot R_T+\nabla\cdot R_N+\nabla\cdot R_H+\nabla\cdot R_S.
\end{split}
\end{equation}
For parameters $\lambda_0\gg 1$, $b>1$ and $0<\beta<1$, define
\[\lambda_q=\left\lceil \lambda_0^{b^q}\right\rceil, \ \ \ q\in \mathbb N\cup \{0\}\]
and let $\delta_q=\lambda_q^{-\beta}$.
Assume $\theta_q$ and $u_q$ are localized to frequency $\sim \lambda_q$. The cancellation $U_{I,q+1}\Theta_{\bar I, q+1}+c_{1,q}A_1$ suggests the scaling $|c_{1,q}|\sim |\theta_{I,q+1}|^2$.
We make the inductive assumptions:
\begin{equation}\label{est-q-1}
\|\nabla^k u_q\|_{C^0}+\|\nabla^k \theta_q\|_{C^0}\lesssim \lambda_q^k\delta_{q-1}^{\frac12}, \ \ k=1, 2, ..., L,
\end{equation}
\begin{equation}\label{est-q-2}
\|\nabla^k (\partial_t+u_q\cdot\nabla)\theta_q\|_{C^0}\lesssim \lambda_q^{k+1}\delta_{q-1}, \ \ k=0, 1, 2, ..., L-1,
\end{equation}
\begin{equation}\label{est-q-3}
\|\nabla^k c_{1,q}\|_{C^0}\lesssim \lambda_q^k\delta_{q}, \ \ k=1, 2, ..., L,
\end{equation}
\begin{equation}\label{est-q-4}
\|\nabla^k (\partial_t+u_q\cdot\nabla)c_{1,q}\|_{C^0}\lesssim \lambda_q^{k+1}\delta_{q-1}^{\frac12}\delta_{q}, \ \ k=0, 1, 2, ..., L-1,
\end{equation}
\begin{equation}\label{est-q-5}
\|\nabla^k R_{q}\|_{C^0}\lesssim \lambda_q^k\delta_{q+1}, \ \ k=1, 2, ..., L,
\end{equation}
\begin{equation}\label{est-q-6}
\|\nabla^k (\partial_t+u_q\cdot\nabla)R_{q}\|_{C^0}\lesssim \lambda_q^{k+1}\delta_{q-1}^{\frac12}\delta_{q+1}, \ \ k=0, 1, 2, ..., L-1.
\end{equation}
The increments $\Theta_{q+1}$ and $U_{q+1}=T[\Theta_{q+1}]$ satisfy
\begin{equation}\label{est-q1-1}
\|\nabla^k\Theta_{q+1}\|_{C^0}+\|\nabla^k U_{q+1}\|_{C^0}\lesssim \lambda_{q+1}^k\delta_q^{\frac12}, \ \ k=0,1,
\end{equation}
\begin{equation}\label{est-q1-2}
\|(\partial_t+u_q\cdot\nabla)\Theta_{q+1}\|_{C^0}+\|(\partial_t+u_q\cdot\nabla) U_{q+1}\|_{C^0}\lesssim \tau_{q}^{-1}\delta_q^{\frac12}
\end{equation}
where the time scale $ \tau_{q}$ is to be determined in the following.

Since we can find $R_T$ and $R_N$ such that
\[R_T=\nabla \Delta^{-1} \mathbb P_{\sim \lambda_{q+1}} [(\partial_t+u_q\cdot\nabla)\Theta_{q+1}]\]
\[R_N=\nabla \Delta^{-1} (U_{q+1}\cdot\nabla \theta_q),\]
it follows from (\ref{est-q-1}), (\ref{est-q1-1}) and (\ref{est-q1-2})
\begin{equation}\notag
\begin{split}
\|R_T\|_{C^0}\lesssim&\ \lambda_{q+1}^{-1} \|(\partial_t+u_q\cdot\nabla)\Theta_{q+1}\|_{C^0}\lesssim  \lambda_{q+1}^{-1}\tau_{q}^{-1}\delta_q^{\frac12},\\
\|R_N\|_{C^0}\lesssim&\ \lambda_{q+1}^{-1} \| U_{q+1}\cdot\nabla \theta_q\|_{C^0}\lesssim  \lambda_{q+1}^{-1}\delta_q^{\frac12}\lambda_q\delta_{q-1}^{\frac12}.
\end{split}
\end{equation}
On the other hand, we have
\begin{equation}\notag
\begin{split}
\|R_H\|_{C^0}=&\left\| \sum_{J\neq \bar I}U_{J,q+1}\Theta_{I, q+1}\right\|_{C^0}\\
\lesssim&\sum_{I}\| \theta_{I, q+1}\|_{C^0}^2\left(m(\|\nabla\xi_I)-m(\nabla\widehat \xi_I)\|_{C^0}+\|\nabla\xi_I-\nabla\widehat \xi_I\|_{C^0} \right)\\
\lesssim&\sum_{I}\| \theta_{I, q+1}\|_{C^0}^2\|\nabla\xi_I-\nabla\widehat \xi_I\|_{C^0} \\
\lesssim&\sum_{I}\| \theta_{I, q+1}\|_{C^0}^2\lambda_q\tau_q\|u_q\|_{C^0} \\
\lesssim&\ \delta_q\lambda_q\delta_{q-1}^{\frac12}\tau_{q}.
\end{split}
\end{equation}
To balance the error $R_T$ and $R_H$, we choose $\tau_{q}=\delta_{q-1}^{-\frac14}\delta_q^{-\frac14}\lambda_q^{-\frac12}\lambda_{q+1}^{-\frac12}$ such that
\[\lambda_{q+1}^{-1}\tau_{q}^{-1}\delta_q^{\frac12}\sim \delta_q\lambda_q\delta_{q-1}^{\frac12}\tau_{q}.\]
In the end, we observe (up to small errors)
\[R_S=c_{2,q+1}A_2\]
for some coefficient $c_{2,q+1}$ with $|c_{2,q+1}|\leq \delta_{q+1}$. We then denote $R_{q+1}=R_T+R_N+R_H$. 
Combining the estimates above gives
\begin{equation}\notag
\begin{split}
\|R_{q+1}\|_{C^0}\lesssim&\ \lambda_{q+1}^{-1}\tau_{q}^{-1}\delta_q^{\frac12}+\lambda_{q+1}^{-1}\delta_q^{\frac12}\lambda_q\delta_{q-1}^{\frac12}\\
\lesssim&\ \lambda_{q+1}^{-\frac12}\lambda_q^{\frac12}\delta_q^{\frac34}\delta_{q-1}^{\frac14}+\lambda_{q+1}^{-1}\delta_q^{\frac12}\lambda_q\delta_{q-1}^{\frac12}\\
\lesssim &\ \lambda_q^{-\frac12b+\frac12-\frac34\beta-\frac{\beta}{4b}}+\lambda_q^{-b+1-\frac12\beta-\frac{\beta}{2b}}.
\end{split}
\end{equation}
To make sure $\|R_{q+1}\|_{C^0}\lesssim \delta_{q+2}$, we require
\begin{equation}\notag
\begin{cases}
-\frac12b+\frac12-\frac34\beta-\frac{\beta}{4b}<-b^2\beta,\\
-b+1-\frac12\beta-\frac{\beta}{2b}<-b^2\beta.
\end{cases}
\end{equation}
Thus we solve, by recalling $b>1$
\begin{equation}\notag
\begin{split}
&b^2\beta-\frac12b+\frac12-\frac34\beta-\frac{\beta}{4b}<0\\
\Longleftrightarrow&\ b\beta(b-1)+\beta(b-1)+\frac{\beta}{4b}(b-1)-\frac12(b-1)<0\\
\Longleftrightarrow&\ b\beta+\beta+\frac{\beta}{4b}-\frac12<0\\
\Longleftrightarrow&\ \beta<\frac{1}{2b+2+\frac{1}{2b}}.
\end{split}
\end{equation}
When $b=1^+$, the inequality above implies $\beta<\frac29$. 
Similarly, the other inequality gives
\begin{equation}\notag
\begin{split}
&b^2\beta-b+1-\frac12\beta-\frac{\beta}{2b}<0\\
\Longleftrightarrow&\ b\beta(b-1)+\beta(b-1)+\frac{\beta}{2b}(b-1)-\frac12(b-1)<0\\
\Longleftrightarrow&\ b\beta+\beta+\frac{\beta}{2b}-1<0\\
\Longleftrightarrow&\ \beta<\frac{1}{b+1+\frac{1}{2b}}
\end{split}
\end{equation}
which indicates $\beta<\frac25$ for $b=1^+$.  The $C^\alpha$ of $\Theta_{q+1}$ requires
\begin{equation}\notag
\|\Theta_{q+1}\|_{C^\alpha}\lesssim \lambda_{q+1}^\alpha \|\Theta_{q+1}\|_{C^0}\lesssim \lambda_{q+1}^\alpha\delta_q^{\frac12}\lesssim \lambda_q^{b\alpha-\frac12\beta}\lesssim 1
\end{equation}
which leads to $\alpha<\frac{\beta}{2b}<\frac19$.

In $d$-dimension, we need to make sure $\|R_{q+1}\|_{C^0}\lesssim \delta_{q+d}$ in order to carry on the iteration, and hence require
\begin{equation}\notag
\begin{cases}
-\frac12b+\frac12-\frac34\beta-\frac{\beta}{4b}<-b^d\beta,\\
-b+1-\frac12\beta-\frac{\beta}{2b}<-b^d\beta.
\end{cases}
\end{equation}
The first inequality is equivalent to
\begin{equation}\notag
\begin{split}
&b^d\beta-b+1-\frac12\beta-\frac{\beta}{2b}<0\\
\Longleftrightarrow&\ b^{d-1}\beta(b-1)+b^{d-2}\beta(b-1)+...+\beta(b-1)+\frac{\beta}{4b}(b-1)-\frac12(b-1)<0\\
\Longleftrightarrow&\ b^{d-1}\beta+b^{d-2}\beta+...+\beta+\frac{\beta}{4b}-\frac12<0\\
\Longleftrightarrow&\ \beta<\frac{1}{2(b^{d-1}+b^{d-2}+...+1+\frac{1}{4b})},
\end{split}
\end{equation}
following which we have $\beta<1/(2d+\frac12)$ for $b>1$. Similarly the second inequality is equivalent to
\begin{equation}\notag
\beta<\frac{1}{b^{d-1}+b^{d-2}+...+1+\frac{1}{2b}}
\end{equation}
and hence $\beta<1/(d+\frac12)$ for $b>1$. Combining the two conditions yields $\beta<\frac{2}{4d+1}$ and hence $\alpha<\frac{\beta}{2b}<\frac{1}{4d+1}$.

\bigskip

\section{Outline of non-uniqueness constructions for forced equation (\ref{ase})}
\label{sec-outline}

In this section we sketch a generic convex integration scheme for forced active scalar equations with even operators. We will explore the flexibility in the convex integration construction due to the presence of an external forcing. Such flexibility was exploited in the previous works \cite{DP, DP2} through the sum-difference formulation of two distinct solutions for SQG. We note an alternating formulation of convex integration techniques was used in 
\cite{BHP2} for forced SQG.

\subsection{Sum-difference system of two solutions} 
\label{sec-two}
Assume $(\theta, u)$ and $(\widetilde \theta, \widetilde u)$ are two distinct solutions of (\ref{ase}). The new variables
\[P=\frac12(\theta+\widetilde\theta), \ \ M=\frac12(\theta-\widetilde\theta)\]
satisfy the system 
\begin{equation}\label{pm1}
\begin{split}
P_t+T[P]\cdot\nabla P+T[M]\cdot\nabla M=&-\nu \Lambda^\gamma P+f,\\
M_t+T[P]\cdot\nabla M+T[M]\cdot\nabla P=&-\nu\Lambda^\gamma M,\\
\nabla\cdot T[P]= 0, \ \ \nabla\cdot T[M]=&\ 0.
\end{split}
\end{equation} 
Allowing forcing in the equation of $M$, we have the flexibility to find a pair $(\theta, u)$ and $(\widetilde\theta, \widetilde u)$ with $\theta-\widetilde \theta\not\equiv 0$, satisfying the relaxed system 
\begin{equation}\label{pm2}
\begin{split}
P_t+T[P]\cdot\nabla P+T[M]\cdot\nabla M=&-\nu \Lambda^\gamma P+f_1,\\
M_t+T[P]\cdot\nabla M+T[M]\cdot\nabla P=&-\nu\Lambda^\gamma M+f_2,\\
\nabla\cdot T[P]= 0, \ \ \nabla\cdot T[M]=&\ 0.
\end{split}
\end{equation} 
We then apply a convex integration scheme to the equation of $M$ with the aim to erase the forcing $f_2$ iteratively.
 

\subsection{The convex integration scheme}

For $f_2\not\equiv 0$, 
we will apply a convex integration scheme to system (\ref{pm2}) with the aim of reducing the forcing $f_2$ in the second equation. We thus consider the approximating system
\begin{equation}\label{pm-q}
\begin{split}
\partial_tP_q+T[P_q]\cdot\nabla P_q+T[M_q]\cdot\nabla M_q=&-\nu \Lambda^\gamma P_q+\nabla\cdot \bar R_q,\\
\partial_tM_q+T[P_q]\cdot\nabla M_q+T[M_q]\cdot\nabla P_q=&-\nu\Lambda^\gamma M_q+\nabla\cdot \widetilde R_q\\
\end{split}
\end{equation} 
inductively. Consistent with notation, we have 
\[\theta_q=P_q+M_q, \ \ \ \widetilde\theta_q=P_q-M_q,\]
\[T[\theta_q]=T[P_q]+T[M_q]=u_q, \ \ \ T[\widetilde\theta_q]=T[P_q]-T[M_q]=\widetilde u_q.\]
Due to the presence of forcing terms in both equations, we have the abundance to find an initial tuple $(P_0, M_0, \bar R_0, \widetilde R_0)$ with $M_0\not\equiv 0$ satisfying (\ref{pm-q}). Starting from this tuple, we construct another solution $(P_1, M_1, \bar R_1, \widetilde R_1)$ of (\ref{pm-q}) with $ \widetilde R_1$ smaller than $\widetilde R_0$ in an appropriate way. Without loss of generality, assume $(P_q, M_q, \bar R_q, \widetilde R_q)$ satisfies (\ref{pm-q}) for an even integer $q$. To take the advantage of the flexibility of having two unknown variables, each stage of the construction consists two steps: from $(P_q, M_q, \bar R_q, \widetilde R_q)$ to $(P_{q+1}, M_{q+1}, \bar R_{q+1}, \widetilde R_{q+1})$ and from $(P_{q+1}, M_{q+1}, \bar R_{q+1}, \widetilde R_{q+1})$ to $(P_{q+2}, M_{q+2}, \bar R_{q+2}, \widetilde R_{q+2})$. In particular, we construct $W_{q+1}$ and $W_{q+2}$ such that
\[M_{q+1}=M_q+W_{q+1}, \ \ \ P_{q+1}=P_q-W_{q+1}\]
and 
\[M_{q+2}=M_{q+1}+W_{q+2}, \ \ \ P_{q+2}=P_{q+1}+W_{q+2}.\]
Consequently we note, for even $q$
\begin{equation}\label{pause-1}
\begin{split}
\theta_{q+1}=&\ P_{q+1}+M_{q+1}=P_{q}+M_{q}=\theta_q,\\
\widetilde\theta_{q+1}=&\ P_{q+1}-M_{q+1}=P_{q}-M_{q}-2W_{q+1}=\widetilde\theta_q-2W_{q+1}
\end{split}
\end{equation}
and 
\begin{equation}\label{pause-2}
\begin{split}
\theta_{q+2}=&\ P_{q+2}+M_{q+2}=P_{q+1}+M_{q+1}+2W_{q+2}=\theta_{q+1}+2W_{q+2},\\
\widetilde\theta_{q+2}=&\ P_{q+2}-M_{q+2}=P_{q+1}-M_{q+1}=\widetilde\theta_{q+1}.
\end{split}
\end{equation}
The ``pause'' reflected in $\theta_{q+1}=\theta_q$ and $\widetilde\theta_{q+2}=\widetilde\theta_{q+1}$ will play a key role to gain better estimates in stress errors.

Since the three tuples  $(P_{q+j}, M_{q+j}, \bar R_{q+j}, \widetilde R_{q+j})$ with $j=0,1,2$ all satisfy (\ref{pm-q}), the stress terms can be expressed as
\begin{equation}\label{R-q1}
\begin{split}
\nabla\cdot \widetilde R_{q+1}=& \left(\partial_t+T[\widetilde \theta_q]\cdot\nabla \right)W_{q+1}+\nu\Lambda^\gamma W_{q+1}
+T[W_{q+1}]\cdot\nabla\widetilde\theta_q\\
&+\left(\nabla\cdot\widetilde R_q-2T[W_{q+1}]\cdot\nabla W_{q+1} \right),
\end{split}
\end{equation}
\begin{equation}\label{R-q2}
\begin{split}
\nabla\cdot \widetilde R_{q+2}=& \left(\partial_t+T[\theta_{q+1}]\cdot\nabla \right)W_{q+2}+\nu\Lambda^\gamma W_{q+2}
+T[W_{q+2}]\cdot\nabla\theta_{q+1}\\
&+\left(\nabla\cdot\widetilde R_{q+1}+2T[W_{q+2}]\cdot\nabla W_{q+2} \right),
\end{split}
\end{equation}
\begin{equation}\label{Rt-q1}
\begin{split}
\nabla\cdot  \bar R_{q+1}=&- \left(\partial_t+T[\widetilde \theta_q]\cdot\nabla \right)W_{q+1}-\nu\Lambda^\gamma W_{q+1}
-T[W_{q+1}]\cdot\nabla\widetilde\theta_q\\
&+\left(\nabla\cdot\bar R_q+2T[W_{q+1}]\cdot\nabla W_{q+1} \right),
\end{split}
\end{equation}
\begin{equation}\label{Rt-q2}
\begin{split}
\nabla\cdot\bar  R_{q+2}=& \left(\partial_t+T[\theta_{q+1}]\cdot\nabla \right)W_{q+2}+\nu\Lambda^\gamma W_{q+2}
+T[W_{q+2}]\cdot\nabla\theta_{q+1}\\
&+\left(\nabla\cdot\bar R_{q+1}+2T[W_{q+2}]\cdot\nabla W_{q+2} \right).
\end{split}
\end{equation}
The forms of the stress terms above provide some insights on the construction of the increments in the two steps $q\to q+1$ and $q+1\to q+2$: \\
(i) $W_{q+j}$ with $j=1,2$ will be designed such that 
\begin{equation}\label{reduce}
\nabla\cdot\widetilde R_q-2T[W_{q+1}]\cdot\nabla W_{q+1} \ \ \mbox{and}\ \ \ \nabla\cdot\widetilde R_{q+1}+2T[W_{q+2}]\cdot\nabla W_{q+2}
\end{equation}
are small;\\
(ii) in view of (\ref{pause-2}), we have $\widetilde\theta_q=\widetilde\theta_{q-1}$ in the iteration process and hence expect to have better estimates for the terms containing $\widetilde\theta_q$ in (\ref{R-q1}); similarly, thanks to $\theta_{q+1}=\theta_q$ in (\ref{pause-1}), better estimates may be achieved for $\widetilde R_{q+2}$ as in (\ref{R-q2});\\
(iii) comparing (\ref{R-q1}) and (\ref{Rt-q1}), the reduced stress error in the process $\widetilde R_q\to \widetilde R_{q+1}$ is gained in the process $R_q\to  R_{q+1}$; while according to (\ref{R-q2}) and (\ref{Rt-q2}), both processes $\widetilde R_{q+1}\to \widetilde R_{q+2}$ and  $R_{q+1}\to R_{q+2}$ have stress error reduced by the same amount. It indicates that this scheme is likely to reduce the forcing in one equation, but not in both equations.

\subsection{Heuristics}
Now we estimate $\widetilde R_{q+1}$ given in (\ref{R-q1}) using the ``pause'' feature $T[\widetilde \theta_q]=T[\widetilde \theta_{q-1}]$ described in item (ii) above. As before, we write 
\[\nabla\cdot\widetilde R_{q+1}=\nabla\cdot R_T+\nabla\cdot R_N+\nabla\cdot R_H+\nabla\cdot R_O+\nabla\cdot R_D\]
with 
\begin{equation}\notag
\begin{split}
\nabla\cdot R_T=&\ \left(\partial_t+T[\widetilde \theta_{q-1}]\cdot\nabla \right)W_{q+1}\\
\nabla\cdot R_N=&\ T[W_{q+1}]\cdot\nabla\widetilde\theta_{q-1}\\
\nabla\cdot R_H=&\ \nabla\cdot \sum_{J\neq \bar I}U_{J,q+1}\Theta_{I, q+1}\\
\nabla\cdot R_D=&\ \nu\Lambda^\gamma W_{q+1}
\end{split}
\end{equation}
and $R_S$ similar as in (\ref{R-q1}). 
Applying (\ref{est-q-1}), (\ref{est-q1-1}) and (\ref{est-q1-2}) gives
\begin{equation}\notag
\begin{split}
\|R_T\|_{C^0}\lesssim&\ \lambda_{q+1}^{-1} \|(\partial_t+\widetilde u_{q-1}\cdot\nabla)\Theta_{q+1}\|_{C^0}\lesssim  \lambda_{q+1}^{-1}\tau_{q}^{-1}\delta_q^{\frac12},\\
\|R_N\|_{C^0}\lesssim&\ \lambda_{q+1}^{-1} \| U_{q+1}\cdot\nabla \widetilde\theta_{q-1}\|_{C^0}\lesssim  \lambda_{q+1}^{-1}\delta_q^{\frac12}\lambda_{q-1}\delta_{q-2}^{\frac12},\\
\|R_D\|_{C^0}\lesssim&\ \lambda_{q+1}^{-1+\gamma} \| W_{q+1}\|_{C^0}\lesssim  \lambda_{q+1}^{-1+\gamma}\delta_q^{\frac12},\\
\|R_H\|_{C^0}\lesssim&\sum_{I}\| \theta_{I, q+1}\|_{C^0}^2\lambda_{q-1}\tau_{q}\|\widetilde u_{q-1}\|_{C^0}
\lesssim \delta_q\lambda_{q-1}\delta_{q-2}^{\frac12}\tau_{q}
\end{split}
\end{equation}
where in the last inequality we used the fact that $\xi_{I}$ is now advected by $\widetilde u_{q-1}$. To balance the estimates of $\|R_T\|_{C^0}$ and $\|R_H\|_{C^0}$, we set
\begin{equation}\notag
\lambda_{q+1}^{-1}\tau_{q}^{-1}\delta_q^{\frac12}=\delta_q\lambda_{q-1}\delta_{q-2}^{\frac12}\tau_{q}
\end{equation}
which implies
\begin{equation}\notag
\tau_{q}=\lambda_{q-1}^{-\frac12}\lambda_{q+1}^{-\frac12}\delta_{q-2}^{-\frac14}\delta_{q}^{-\frac14}.
\end{equation}
In the end, to obtain 
\[ \|R_{q+1}\|_{C^0}\lesssim \delta_{q+2}\]
we need to require 
\begin{equation}\notag
\begin{cases}
\lambda_{q+1}^{-1}\tau_{q}^{-1}\delta_q^{\frac12}\lesssim \delta_{q+2},\\
\lambda_{q+1}^{-1}\delta_q^{\frac12}\lambda_{q-1}\delta_{q-2}^{\frac12} \lesssim \delta_{q+2},\\
\lambda_{q+1}^{-1+\gamma}\delta_q^{\frac12} \lesssim \delta_{q+2}
\end{cases}
\end{equation}
which are valid provided
\begin{equation}\notag
\begin{split}
&b^2\beta-\frac34\beta-\frac{\beta}{4b^2}+\frac{1}{2b}-\frac{b}{2}<0\\
\Longleftrightarrow&\ b\beta(b-1)+\beta(b-1)+\frac{\beta}{4b^2}(b^2-1)-\frac{1}{2b}(b^2-1)<0\\
\Longleftrightarrow&\ b\beta+\beta+\frac{\beta}{4b^2}(b+1)-\frac{1}{2b}(b+1)<0\\
\Longleftrightarrow&\ \beta<\frac{\frac{1}{2b}(b+1)}{b+1+\frac{b+1}{4b^2}}
\end{split}
\end{equation}
and 
\begin{equation}\notag
b(-1+\gamma)-\frac12\beta+b^2\beta<0 \Longleftrightarrow \beta<\frac{2b(1-\gamma)}{2b^2-1}.
\end{equation}
So for $b=1^+$, we have $\beta<\frac25$, $\alpha<\frac{\beta}{2b}<\frac15$, and $\gamma<1-\alpha$.

In general for $d$-dimension, we need to impose 
\[ \|R_{q+1}\|_{C^0}\lesssim \delta_{q+d}\]
and hence 
\begin{equation}\notag
\begin{cases}
\lambda_{q+1}^{-1}\tau_{q}^{-1}\delta_q^{\frac12}\lesssim \delta_{q+d},\\
\lambda_{q+1}^{-1}\delta_q^{\frac12}\lambda_{q-1}\delta_{q-2}^{\frac12} \lesssim \delta_{q+d},\\
\lambda_{q+1}^{-1+\gamma}\delta_q^{\frac12} \lesssim \delta_{q+d}.
\end{cases}
\end{equation}
Thus we have
\begin{equation}\notag
\begin{split}
&b^d\beta-\frac34\beta-\frac{\beta}{4b^2}+\frac{1}{2b}-\frac{b}{2}<0\\
\Longleftrightarrow&\ b^{d-1}\beta(b-1)+...+\beta(b-1)+\frac{\beta}{4b^2}(b^2-1)-\frac{1}{2b}(b^2-1)<0\\
\Longleftrightarrow&\ b^{d-1}\beta+...+b\beta+\beta+\frac{\beta}{4b^2}(b+1)-\frac{1}{2b}(b+1)<0\\
\Longleftrightarrow&\ \beta<\frac{\frac{1}{2b}(b+1)}{b^{d-1}+...+b+1+\frac{b+1}{4b^2}}
\end{split}
\end{equation}
and 
\begin{equation}\notag
b(-1+\gamma)-\frac12\beta+b^d\beta<0 \Longleftrightarrow \beta<\frac{2b(1-\gamma)}{2b^d-1}.
\end{equation}
For $b=1^+$, $\beta<\frac{1}{d+\frac12}$, $\alpha<\frac{\beta}{2b}<\frac{1}{2d+1}$, and $\gamma<1-\alpha$.

\medskip
\subsection{Key idea to reduce the stress error}
\label{sec-couple}
We discuss how to achieve the error reduction in (\ref{reduce}). The assumptions on the Fourier symbol $m$ of the operator $T$ imply that there are two linearly independent vectors in the image of the even part of $m$
\begin{equation}\label{two-vec}
A_1=m(\xi^{(1)})+m(-\xi^{(1)}), \ \ \ A_2=m(\xi^{(2)})+m(-\xi^{(2)})
\end{equation}
with $\xi^{(1)}, \xi^{(2)}\in \mathbb Z^2$. 

Consider the increment ansatz 
\begin{equation}\notag
\begin{split}
W_{I,q+1}=&\ \mathbb P_{I,\lambda_{q+1}}\left(a_{I,q+1}(x,t)e^{i\lambda_{q+1}\xi_I}\right),  \\ 
W_{q+1}=&\sum_{I\in\Omega} W_{I,q+1}.
\end{split}
\end{equation}
By the Microlocal Lemma \ref{le-mic} and zero degree of homogeneity of $m$, the drift term takes the form
\begin{equation}\notag
T[W_{I,q+1}]=m(\nabla \xi_I) W_{I, q+1}+\delta T[W_{I,q+1}]
\end{equation}
with a small error term $\delta T[W_{I,q+1}]$. A straightforward computation shows that
\begin{equation}\notag
\begin{split}
T[W_{q+1}] W_{q+1}=&\frac12\sum_{I\in\Omega} \left(T[W_{I,q+1}] W_{\bar I, q+1}+T[W_{\bar I,q+1}] W_{I, q+1}\right)\\
&+\sum_{J\neq \bar I} T[W_{I, q+1}] W_{J, q+1}\\
=&\ \frac12\sum_{I\in\Omega}|a_{I,q+1}|^2\left(m(\nabla\xi_I)+ m(-\nabla\xi_I)\right)+\mbox{error}\\
&+\sum_{J\neq \bar I} T[W_{I, q+1}] W_{J, q+1}.
\end{split}
\end{equation}
Since $m$ is not odd, the leading order (low frequency) term is non zero, 
\[|a_{I,q+1}|^2\left(m(\nabla\xi_I)+ m(-\nabla\xi_I)\right)\neq 0.\]
The goal is to construct coefficient functions $a_{I, q+1}$ and phase functions $\xi_I$ such that 
\begin{equation}\label{reduce-Rq}
\nabla\cdot \left(\widetilde R_q-\sum_{I\in\Omega}|a_{I,q+1}|^2\left(m(\nabla\xi_I)+ m(-\nabla\xi_I)\right)\right) 
\end{equation}
is small. As $A_1$ and $A_2$ defined in (\ref{two-vec}) span $\mathbb R^2$, we expect to choose $\xi_I$ to guarantee 
\[m(\nabla\xi_I)+ m(-\nabla\xi_I)\approx A_1 \ \ \mbox{or} \ \ m(\nabla\xi_I)+ m(-\nabla\xi_I)\approx A_2. \]
Then with an appropriate choice of $a_{I,q+1}$ we hope the principal part of $\widetilde R_q$ can be canceled through (\ref{reduce-Rq}). The reduction of $\widetilde R_{q+1}$ is achieved analogously. 

To ensure other terms in the stress fields given in (\ref{R-q1}) and (\ref{R-q2}) can be controlled appropriately, we need several technical tools which will be provided in Section \ref{sec-pre}.


\bigskip

\section{Construction of the increment}
\label{sec-pre}



We describe the construction of the highly oscillatory correction (increment) in this section. We start with some technical preparations. 

\subsection{Microlocal Lemma}
The drift operator $T$ is a nonlocal differential operator. When acting $T$ on plane waves, we need the following lemma from \cite{IV} to extract the leading order term. 

\begin{Lemma}[Microlocal Lemma] \label{le-mic}
Let $K:\mathbb R^2\to\mathbb C$ be a Schwartz function and 
\[T[\Theta](x)=\int_{\mathbb R^2}\Theta(x-h) K(h)\, dh\]
for $\Theta:\mathbb T^2\to\mathbb C$. For any $\Theta=e^{i\lambda\xi(x)}\theta(x)$ with $\lambda\in \mathbb Z$ and smooth functions $\xi:\mathbb T^2\to\mathbb R$ and $\theta:\mathbb T^2\to\mathbb C$, we have
\begin{equation}\notag
T[\Theta](x)=e^{i\lambda\xi(x)}\left(\theta(x)\widehat K(\lambda\nabla\xi)+\delta[T\Theta](x) \right)
\end{equation}
with the error term given by
\begin{equation}\notag
\begin{split}
\delta[T\Theta](x) =& \int_0^1\frac{d}{dr} \int_{\mathbb R^2}e^{-i\lambda\nabla \xi\cdot h}e^{iZ(r,x,h)}\theta(x-rh)K(h) \, dr\, dh\\
Z(r,x,h)=&\ r\lambda\int_0^1h^jh^l\partial_j\partial_l \xi(x-sh)(1-s)\, ds.
\end{split}
\end{equation}
\end{Lemma}



\medskip

\subsection{Mollification}
As is standard in convex integration method, to avoid loss of derivative, we need to regularize the solution $(P_q, M_q, \bar R_q, \widetilde R_q)$ before adding increments to produce $(P_{q+1}, M_{q+1}, \bar R_{q+1}, \widetilde R_{q+1})$. We first regularize $P_q, M_q, T[P_q]$ and $T[M_q]$. Fix some $L\geq 1$. Choose $\mu_q=\lambda_{q+1}^{\frac{1}{L}}\lambda_q^{1-\frac{1}{L}}$. Denote $\mathbb P_{\leq \mu_q}$ by the Littlewood-Paley projection onto frequency $\leq \mu_q$. Define 
\begin{equation}\notag
\begin{split}
P_{\epsilon, q}=&\ \mathbb P^2_{\leq \mu_q} P_q, \ \ M_{\epsilon, q}=\mathbb P^2_{\leq \mu_q} M_q, \\ 
T[P_{\epsilon, q}]=&\ \mathbb P^2_{\leq \mu_q} T[P_q], \ \ T[M_{\epsilon, q}]=\mathbb P^2_{\leq \mu_q} T[M_q].
\end{split}
\end{equation}

\begin{Lemma}\label{le-molli}
The estimates  
\begin{equation}\notag
\begin{split}
\|P_{q}-P_{\epsilon, q}\|_{C^0}\lesssim&\ \mu_q^{-j}\|\nabla^j P_{q}\|_{C^0},\\
\|M_{q}-M_{\epsilon, q}\|_{C^0}\lesssim&\ \mu_q^{-j}\|\nabla^j M_{q}\|_{C^0},\\
\|T[P_{q}]-T[P_{\epsilon, q}]\|_{C^0}\lesssim&\ \mu_q^{-j}\|\nabla^j T[P_{q}]\|_{C^0},\\
\|T[M_{q}]-T[M_{\epsilon, q}]\|_{C^0}\lesssim&\ \mu_q^{-j}\|\nabla^j T[M_{q}]\|_{C^0}
\end{split}
\end{equation}
hold for $0\leq j\leq L$.
\end{Lemma}

The regularization of $\widetilde R_q$ is slightly more involved. We first define the coarse scale flow $ \Phi_q(x, s; t): \mathbb T^2\times\mathbb R\times\mathbb R\to  \mathbb T^2\times\mathbb R$ to be the solution of 
\begin{equation}\notag
\begin{cases}
\partial_s \Phi_q(x,s;t)=\widetilde u_q( \Phi_q(x,s;t), s)\\
\Phi_q(x,t;t)=x
\end{cases}
\end{equation}
with $\widetilde u_q=T[M_q]$. 
Let $\eta$ be a standard mollifier in space with $\supp \eta\subset B(0,1)$ 
and $\rho$ be a standard mollifier in time with $\supp \rho\subset (-1,1)$. Denote 
\[\eta_\delta(x)=\delta^{-d}\eta(\delta^{-1}x), \ \ \rho_\delta(s)=\delta^{-1}\rho(\delta^{-1}s).\] 
Without loss of generality, assume $\widetilde R_q$ can be decomposed as
\[\widetilde R_q=c_{1,q}A_1+R^*_q.\]
For the spatial scale $\e_x$ and time scale $\e_t$ to be determined later, we define the regularized stress and component
\begin{equation}\notag
\begin{split}
R^*_{q, \e_x}=&\ \eta_{\e_x}* R^*_{q},\\
R^*_{q, \e}(x,t)=&\int_{\mathbb R} R^*_{q, \e_x}(\Phi_q(x,t+s;t), t+s) \rho_{\e_t} (s)\, ds,\\
 c_{1, q, \e_x}=&\ \eta_{\e_x}*  c_{1,q},\\
 c_{1, q, \e}(x,t)=&\int_{\mathbb R} c_{1, q, \e_x}(\Phi_q(x,t+s;t), t+s) \rho_{\e_t} (s)\, ds,
\end{split}
\end{equation}
and 
\[\widetilde R_{q, \e}=c_{1,q, \e}A_1+R^*_{q,\e}. \]
The purpose of such regularization is to have better estimates on the advective derivatives, since the advective derivative commutes with the flow map $\Phi$. See \cite{Is0} for more details. 

Choose the length and time scales
\begin{equation}\label{molli-scale}
\e_x=(\lambda_{q+1}^{-1}\lambda_q)^{\frac1{L}}\lambda_q^{-1}, \ \ \e_t=\lambda_{q+1}^{-1}\delta_q^{-\frac12}.
\end{equation}

\begin{Lemma}\label{le-molli-est}
The regularized stress field satisfies the following estimates
\begin{equation}\notag
\begin{split}
\|(c_{1, q, \e}-c_{1,q})A_1\|_{C^0}+\|R^*_{q,\e}- R^*_{q}\|_{C^0}\lesssim \delta_{q-1}^{\frac12}\delta_{q}^{\frac12}\lambda_{q+1}^{-1}\lambda_q,\\
\|\nabla^k c_{1, q, \e}\|_{C^0}\lesssim \lambda_q^k\delta_q(\lambda_{q+1}\lambda_q^{-1})^{\frac{(k+1-L)_+}{L}},\\
\|\nabla^k R^*_{ q, \e}\|_{C^0}\lesssim \lambda_q^k\delta_{q+1}(\lambda_{q+1}\lambda_q^{-1})^{\frac{(k+1-L)_+}{L}},\\
\|\nabla^k \frac{D}{Dt}c_{1, q, \e}\|_{C^0}\lesssim \lambda_q^k\delta_q\lambda_q\delta_{q-1}^{\frac12}(\lambda_{q+1}\lambda_q^{-1})^{\frac{(k+1-L)_+}{L}},\\
\|\nabla^k \frac{D}{Dt}R^*_{q, \e}\|_{C^0}\lesssim \lambda_q^k\delta_{q+1}\lambda_q\delta_{q-1}^{\frac12}(\lambda_{q+1}\lambda_q^{-1})^{\frac{(k+1-L)_+}{L}},\\
\|\nabla^k (\frac{D}{Dt})^2c_{1, q, \e}\|_{C^0}\lesssim \lambda_q^k\delta_q\lambda_q\delta_{q-1}^{\frac12}\lambda_{q+1}\delta_q^{\frac12}(\lambda_{q+1}\lambda_q^{-1})^{\frac{(k+1-L)_+}{L}},\\
\|\nabla^k (\frac{D}{Dt})^2R^*_{q, \e}\|_{C^0}\lesssim \lambda_q^k\delta_{q+1}\lambda_q\delta_{q-1}^{\frac12}\lambda_{q+1}\delta_q^{\frac12}(\lambda_{q+1}\lambda_q^{-1})^{\frac{(k+1-L)_+}{L}}.
\end{split}
\end{equation}
\end{Lemma}

\medskip

\subsection{Lifting function}
To achieve the cancellation in (\ref{reduce-Rq}), we further decompose $ R^*_{q,\e}$ as 
\begin{equation}\label{c2-q1}
R^*_{q,\e}=c_{1,q, r}A_1+c_{2,q+1}A_2
\end{equation}
and hence
\begin{equation}\notag
\widetilde R_{q,\e}=\left(c_{1, q, \e}+c_{1,q, r}\right)A_1+c_{2,q+1}A_2.
\end{equation}
Then we require 
\begin{equation}\notag
\sum_{I}|a_{I,q+1}|^2 A_1=e_q(t)A_1+\left(c_{1, q, \e}+c_{1,q, r}\right)A_1
\end{equation}
for some function $e_q(t)$ such that 
\[e_q(t)+\left(c_{1, q, \e}+c_{1,q, r}\right)>0, \ \ \ e_q(t)>2\left(c_{1, q, \e}+c_{1,q, r}\right).\] 
The function $e_q$ is referred to be the lifting function.

\medskip

\subsection{Time cutoff}
To optimize the control of main error terms, we choose the life span $\tau_{q}$ of the increment $W_{I,q+1}$ to be an appropriate short time interval. The size of $\tau_{q}$ will be specified later. Let $\phi$ satisfy the condition of partition of unity in time
\begin{equation}\notag
\sum_{n\in \mathbb Z} \phi^2(t-n)=1.
\end{equation}
Define 
\begin{equation}\notag
\phi_k(t)=\phi\left( \frac{t-k\tau_{q}}{\tau_{q}}\right)
\end{equation}
and consider the amplitude function 
\begin{equation}\notag
a_{I,q+1}=e^{\frac12}(t)\phi_k(t)(1+(c_{1,q,\e}+c_{1,q,r})e^{-1}(t))^{\frac12}, \ \ I=(k,\pm).
\end{equation}
The choice of $\phi_k$ indicates that the amplitude function $a_{I, q+1}$ with $I=(k,f)$ has support $[k\tau_{q}-\frac23\tau_{q}, k\tau_{q}+\frac23\tau_{q}]$ in time.  

\medskip

\subsection{Phase functions}
Denote $\widetilde u_{q, \e}=T[P_{q,\e}]-T[M_{q,\e}]$.
We identify the phase functions $\xi_I$ to be solutions of the transport equation
\begin{equation}\label{transport}
\begin{split}
(\partial_t+\widetilde u_{q, \e}\cdot\nabla)\xi_I=&\ 0\\
\xi_I(k\tau_{q}, x)=&\ \widehat\xi_I(x), 
\end{split}
\end{equation}
where the initial data for the phase function is given by
\begin{equation}\notag
\widehat\xi_I(x)=\widehat\xi_{(k,\pm)}(x)=\pm10^{[k]}\xi^{(1)}\cdot x
\end{equation}
with 
\begin{equation}\notag
[k]=
\begin{cases}
0, \ \  \mbox{if } \ \ k \ \ \mbox{is even},\\
1, \ \  \mbox{if } \ \ k \ \ \mbox{is odd}.
\end{cases}
\end{equation}
In the iteration step of $q+1\to q+2$, the initial data is chosen as 
\begin{equation}\notag
\widehat\xi_I(x)=\widehat\xi_{(k,\pm)}(x)=\pm10^{[k]}\xi^{(2)}\cdot x.
\end{equation}
The time scale $\tau_q$ will be chosen such 
\[|\nabla\xi_I-\nabla\widehat \xi_I|\leq \frac14 |\nabla\widehat \xi_I|.\]

\medskip

\subsection{Construction of the increment $W_{q+1}$} \label{subsec-construct}
We are ready to introduce the increment $W_{q+1}$,
\begin{equation}\label{increment-q}
W_{q+1}=\sum_{I\in \Omega} W_{I, q+1}, \ \ \ W_{I,q+1}=\mathbb P_{\approx \lambda_{q+1}} [a_{I,q+1} e^{i\lambda_{q+1}\xi_I}]
\end{equation}
with the amplitude function
\begin{equation}\label{amplitude}
a_{I,q+1}(x,t)=e^{\frac12}(t)\phi_k(t) \left(1+(c_{1,q,\e}+c_{1,q,r})(x,t)e^{-1}(t)\right)^{\frac12} , \ \ I=(k,\pm).
\end{equation}
The projection operator $\mathbb P_{\approx \lambda_{q+1}}$ is defined through its Fourier multiplier as
\begin{equation}\notag
\widehat {\mathbb P_{\approx \lambda_{q+1}} f}(\xi)=\eta_{\lambda_{q+1}}(\xi)\hat f(\xi).
\end{equation}
We specify the multiplier $\eta_{\lambda_{q+1}}$ in the following. Let $\eta$ be a smooth bump function in Fourier space with frequency support on $B_{|\xi^{(1)}|/2}(\xi^{(1)})$, satisfying
\begin{equation}\notag
\eta(\xi)=1, \ \ \mbox{if} \ \ |\xi-\xi^{(1)}|\leq \frac14|\xi^{(1)}|.
\end{equation}
Then we define 
\begin{equation}\notag
\eta_{\lambda_{q+1}}(\xi)=\eta(\pm 10^{-[k]}\lambda_{q+1}^{-1}\xi).
\end{equation}
Applying the Microlocal Lemma \ref{le-mic} gives
\begin{equation}\notag
\begin{split}
W_{I,q+1}=&\ e^{i\lambda_{q+1}\xi_I}(a_{I,q+1}+\delta a_{I,q+1}),\\
T[W_{I,q+1}]=&\ e^{i\lambda_{q+1}\xi_I}(u_{I,q+1}+\delta u_{I,q+1})
\end{split}
\end{equation}
with $u_{I,q+1}=m(\lambda_{q+1}\nabla\xi_I)a_{I,q+1}=m(\nabla\xi_I)a_{I,q+1}$ since $m$ is homogeneous of order 0, where $\delta a_{I,q+1}$ and $\delta u_{I,q+1}$ are error terms.

\medskip

\subsection{Main iteration process}

For the tuple $(P_q, M_q, \bar R_q, \widetilde R_q)$ satisfying (\ref{pm-q}), we make the following inductive assumptions. 
Assume $ \widetilde R_q$ can be written as 
\begin{equation}\label{Rq-decomp}
\widetilde R_q=
\begin{cases}
c_{1,q}A_1+R^*_{q}, \ \ \mbox{if} \ \ q \ \ \mbox{is even} \\
 c_{2,q}A_2+R^*_{q}, \ \ \mbox{if} \ \ q \ \ \mbox{is odd}.
\end{cases}
\end{equation}
Assume the estimates below hold:
\begin{equation}\label{induct-q-1}
\begin{split}
&\|\nabla^k P_q\|_{C^0}+\|\nabla^k M_q\|_{C^0}\\
&+\|\nabla^k T[P_q]\|_{C^0}+\|\nabla^k T[M_q]\|_{C^0}\lesssim \lambda_q^{k}\delta_{q-1}^{\frac12}, \ \ k=1, ..., L
\end{split}
\end{equation}
\begin{equation}\label{induct-q-2}
\begin{split}
&\|\nabla^k (\partial_t+T[P_q]\cdot\nabla)T[M_q]\|_{C^0}\\
&+\|\nabla^k (\partial_t+T[M_q]\cdot\nabla)T[M_q]\|_{C^0}\lesssim \lambda_q^{k+1}\delta_{q-1}, \ \ k=0, 1, ..., L-1
\end{split}
\end{equation}
\begin{equation}\label{induct-q-3}
\|\nabla^k c_{1,q}\|_{C^0}\lesssim \lambda_q^{k}\delta_{q}, \ \ k=0, 1, ..., L
\end{equation}
\begin{equation}\label{induct-q-4}
\begin{split}
&\|\nabla^k (\partial_t+T[P_q]\cdot\nabla)c_{1,q}\|_{C^0}\\
&+\|\nabla^k (\partial_t+T[M_q]\cdot\nabla)c_{1,q}\|_{C^0}\lesssim \lambda_q^{k+1}\delta_{q-1}^{\frac12}\delta_q, \ \ k=0, 1, ..., L-1
\end{split}
\end{equation}
\begin{equation}\label{induct-q-5}
\|\nabla^k R^*_{q}\|_{C^0}\lesssim \lambda_q^{k}\delta_{q+1}, \ \ k=0, 1, ..., L
\end{equation}
\begin{equation}\label{induct-q-6}
\begin{split}
&\|\nabla^k (\partial_t+T[P_q]\cdot\nabla)R^*_{q}\|_{C^0}\\
&+\|\nabla^k (\partial_t+T[M_q]\cdot\nabla)R^*_{q}\|_{C^0}\lesssim \lambda_q^{k+1}\delta_{q-1}^{\frac12}\delta_{q+1}, \ \ k=0, 1, ..., L-1.
\end{split}
\end{equation}

In view of the crucial feature (\ref{pause-1})-(\ref{pause-2}) of the two-step scheme, we have for even $q$ and $k=0, 1, ..., L$
\begin{equation}\label{induct-q-7}
\begin{cases}
\widetilde \theta_q=\widetilde \theta_{q-1}, \ \ \|\nabla^k \widetilde \theta_q\|_{C^0}\lesssim \lambda_{q-1}^{k}\delta_{q-2}^{\frac12}, \\ 
\theta_{q+1}=\theta_{q}, \ \ \|\nabla^k \theta_q\|_{C^0}\lesssim \lambda_{q}^{k}\delta_{q-1}^{\frac12}.
\end{cases}
\end{equation}

\begin{Proposition}[Main Iteration] \label{prop-main}
Let $L\geq 2$ and $K,C\geq 4$. Assume $(P_q, M_q, \bar R_q, \widetilde R_q)$ satisfies (\ref{pm-q}) and (\ref{Rq-decomp})-(\ref{induct-q-7}). Let $I_t\subset\mathbb R$ be a nonempty closed interval such that 
\begin{equation}\label{time-support}
\supp R^*_q\cup \supp c_{1,q}\subset I_t\times \mathbb T^2.
\end{equation} 
Let $e$ be a function of time satisfying
\begin{equation}\label{e-support}
e(t)\geq K \delta_q \ \ \forall \ \ t\in I_t\pm \hat \tau_q
\end{equation}
with the natural time scale $\hat \tau_q=\lambda_q^{-1}\delta_{q-1}^{-\frac12}$, and 
\begin{equation}\label{est-lift}
\|\frac{d^r}{dt^r}e^{\frac12}(t)\|_{C^0}\leq C (\lambda_q\delta_{q-1}^{\frac12})^r\delta_q^{\frac12}, \ \ 0\leq r\leq 2.
\end{equation}
There exists another tuple $(P_{q+1}, M_{q+1}, \bar R_{q+1}, \widetilde R_{q+1})$ satisfying (\ref{pm-q}) with $q$ replaced by $q+1$ in the form 
\begin{equation}\notag
\begin{cases}
M_{q+1}=M_q+W_{q+1}, \ \ P_{q+1}=P_q-W_{q+1}, \ \ \mbox{if} \ \ q \ \ \mbox{is even},\\
M_{q+1}=M_q+W_{q+1}, \ \ P_{q+1}=P_q+W_{q+1}, \ \ \mbox{if} \ \ q \ \ \mbox{is odd}
\end{cases}
\end{equation}
and $W_{q+1}=\nabla\cdot \widetilde W_{q+1}$.
Moreover, $\widetilde R_{q+1}$ can be written as 
\begin{equation}\label{Rq1-decomp}
\widetilde R_{q+1}=
\begin{cases}
c_{2,q+1}A_2+R^*_{q+1}, \ \ \mbox{if} \ \ q \ \ \mbox{is even}, \\
 c_{1,q+1}A_1+R^*_{q+1}, \ \ \mbox{if} \ \ q \ \ \mbox{is odd}
\end{cases}
\end{equation}
with 
\begin{equation}\label{time-support-Rq1}
(\supp R^*_{q+1}\cup \supp c_{j,q+1})\subset \supp e\times \mathbb T^2, \ \ j=1,2.
\end{equation} 
The estimates (\ref{induct-q-1})-(\ref{induct-q-6}) are satisfied with $q$ replaced by $q+1$ and in particular $c_{1,q}$ replaced by $c_{2, q+1}$. 
The correction $W_{q+1}$ satisfies the estimates
\begin{equation}\label{est-correct-1}
\|\nabla^k W_{q+1}\|_{C^0}+\|\nabla^k T[W_{q+1}]\|_{C^0}\lesssim \lambda_{q+1}^k\delta_q^{\frac12}, \ \ k=0,1,
\end{equation}
\begin{equation}\label{est-correct-2}
\begin{split}
&\|(\partial_t+T[P_{\e,q}]\cdot\nabla)W_{q+1}\|_{C^0}
+\|(\partial_t+T[P_{\e,q}]\cdot\nabla) T[W_{q+1}]\|_{C^0}\\
&+\|(\partial_t+T[M_{\e,q}]\cdot\nabla)W_{q+1}\|_{C^0}
+\|(\partial_t+T[M_{\e,q}]\cdot\nabla) T[W_{q+1}]\|_{C^0}\\
\lesssim&\ \tau_{q}^{-1}\delta_q^{\frac12},
\end{split}
\end{equation}
\begin{equation}\label{est-correct-3}
\|\nabla^k\widetilde W_{q+1}\|_{C^0}\lesssim \lambda_{q+1}^{k-1}\delta_q^{\frac12}, \ \ k=0,1,
\end{equation}
\begin{equation}\label{est-correct-4}
\begin{split}
&\|(\partial_t+T[P_{\e,q}]\cdot\nabla)\widetilde W_{q+1}\|_{C^0}\\
&+\|(\partial_t+T[M_{\e,q}]\cdot\nabla)\widetilde W_{q+1}\|_{C^0}\lesssim \tau_{q}^{-1}\lambda_{q+1}^{-1}\delta_q^{\frac12}.
\end{split}
\end{equation}
\end{Proposition}

\bigskip

\section{Proof of the main iteration argument}
\label{sec-iteration}

\medskip

\subsection{Basic estimates of increments}

\begin{Lemma}\label{le-est-adv}
Let $L\geq 2$ be the integer in Proposition \ref{prop-main}. The regularizations $T[P_{\e, q}]$ and $T[M_{\e, q}]$ satisfy
\begin{equation}\notag
\begin{split}
&\|\nabla^k T[P_{\e,q}]\|_{C^0}+\|\nabla^k T[M_{\e,q}]\|_{C^0}\\
\lesssim&\ \lambda_q^k\delta_{q-1}^{\frac12} (\lambda_{q+1}\lambda_q^{-1})^{\frac{(k-L)_+}{L}}, \ \ k\geq 1,\\
&\|\nabla^k (\partial_t+T[P_{\e,q}]\cdot\nabla)T[M_{\e,q}]\|_{C^0}\\
&+\|\nabla^k (\partial_t+T[M_{\e,q}]\cdot\nabla)T[M_{\e,q}]\|_{C^0}\\
\lesssim&\ \lambda_q^{k+1}\delta_{q-1} (\lambda_{q+1}\lambda_q^{-1})^{\frac{(k+1-L)_+}{L}}, \ \ k\geq 0.
\end{split}
\end{equation}
\end{Lemma}
See \cite{IV} (Lemma 7.1) for a proof.

\begin{Lemma}\label{le-commu}
Let $h$ be a kernel function satisfying 
\begin{equation}\notag
\||x|^a|\nabla^b h|(x)\|_{L^1(\mathbb R^2)}\leq \lambda^{b-a}, \ \ \lambda\geq \lambda_{q+1}, \ \ 0\leq a\leq b\leq N.
\end{equation} 
Denote 
\[\frac{D_{P,q,\e}}{Dt}=\partial_t+T[P_{\e, q}]\cdot\nabla, \ \ \frac{D_{M,q,\e}}{Dt}=\partial_t+T[M_{\e, q}]\cdot\nabla.\]
For the convolution operator 
\begin{equation}\notag
Qf(x)=\int_{\mathbb R^2} f(x-y)h(y)\, dy,
\end{equation}
the commutators $\left[ \frac{D_{P,q,\e}}{Dt}, Q\right]$ and $\left[ \frac{D_{M,q,\e}}{Dt}, Q\right]$
are bounded operators on $C^0(\mathbb T^2\times \mathbb R)$ and satisfy
\begin{equation}\notag
\left\| \nabla^k[\frac{D_{P,q,\e}}{Dt}, Q]\right\|+\left\| \nabla^k[\frac{D_{M,q,\e}}{Dt}, Q]\right\|\lesssim \lambda_q\delta_{q-1}^{\frac12}\lambda^k, \ \ 0\leq k\leq N-1.
\end{equation}
\end{Lemma}

\begin{Lemma}\label{le-phase}
Let $L\geq 2$ be the integer in Proposition \ref{prop-main}. 
For $\frac{D_q}{Dt}\in \{\frac{D_{P,q}}{Dt}, \frac{D_{M,q}}{Dt} \}$, define
\begin{equation}\notag
D_q^{(k,r)}=\nabla^{k_1}(\frac{D_q}{Dt})^{r_1}\nabla^{k_2}(\frac{D_q}{Dt})^{r_2}\nabla^{k_3}, \ k_1+k_2+k_3=k, \ r_1+r_2=r.
\end{equation}
The phase function $\xi_I$ satisfies 
\begin{equation}\notag
\begin{split}
&\|\nabla^k(\frac{D_q}{Dt})^r\nabla\xi_I\|_{C^0}+\|D_q^{(k,r)}\nabla\xi_I\|_{C^0}\\
\lesssim&\ \lambda_q^k(\lambda_q\delta_{q-1}^{\frac12})^r (\lambda_{q+1}\lambda_q^{-1})^{\frac{(k+1+(r-1)_+-L)_+}{L}}, \ \ k\geq 1, \ r=0,1,2.
\end{split}
\end{equation}
Moreover, we have
\begin{equation}\notag
|\nabla\xi_I(\Phi_q(x,s;t)) -\nabla\widehat\xi_{I}(x)|\leq C \tau_{q}\lambda_q\delta_{q-1}^{\frac12}, \ \ |s|\leq \tau_{q}.
\end{equation}
\end{Lemma}

\begin{Lemma}\label{le-amp1}
Let $L\geq 2$. 
The principal part of the amplitude function satisfies the estimate 
\begin{equation}\notag
\|D_{q'}^{(k,r)}a_{I,q+1}\|_{C^0}+\|D_{q'}^{(k,r)}u_{I,q+1}\|_{C^0}\lesssim \lambda_{q'}^k\delta_q^{\frac12} \tau_{q'}^{-r} (\lambda_{q'+1}\lambda_{q'}^{-1})^{\frac{(k+1-L)_+}{L}}
\end{equation}
for $q'\geq 0$,  $k\geq 0$ and $r=0,1,2$.
\end{Lemma}

\begin{Lemma}\label{le-amp2}
Let $L\geq 2$. The amplitude error terms $\delta a_{I,q+1}$ and $\delta u_{I,q+1}$ satisfy 
\begin{equation}\notag
\|D_{q'}^{(k,r)}\delta a_{I,q+1}\|_{C^0}+\|D_{q'}^{(k,r)}\delta u_{I,q+1}\|_{C^0}\lesssim (\lambda_{q'+1}^{-1}\lambda_{q'})\lambda_q^k\delta_q^{\frac12} \tau_{q'}^{-r} (\lambda_{q+1}\lambda_q^{-1})^{\frac{(k+2-L)_+}{L}}
\end{equation}
for $q'\geq 0$, $k\geq 0$ and $r=0,1,2$.
\end{Lemma}

\begin{Lemma}\label{le-correct}
The estimates for the corrections $W_{q+1}$ and $T[W_{q+1}]$
\begin{equation}\notag
\|D_{q'}^{(k,r)}W_{I,q+1}\|_{C^0}+\|D_{q'}^{(k,r)}T[W_{I,q+1}]\|_{C^0}\lesssim \lambda_{q'+1}^k\tau_{q'}^{-r}\delta_q^{\frac12}
\end{equation}
hold for $q'\geq 0$, $k\geq 0$ and $r=0,1,2$.
\end{Lemma}
The lemmas above can be proved analogously as in \cite{IV}.

\medskip

\subsection{Proof of Proposition \ref{prop-main}}
We only prove the statements for even $q$; the statements for odd $q$ can be established by minor modifications of the proof. Let $W_{q+1}$ be the correction term constructed in Subsection \ref{subsec-construct} and define 
\[M_{q+1}=M_q+W_{q+1}, \ \ P_{q+1}=P_q-W_{q+1}.\]
For $\bar R_{q+1}$ and $\widetilde R_{q+1}$ defined respectively through (\ref{Rt-q1}) and (\ref{R-q1}), the tuple 
\[(P_{q+1}, M_{q+1}, \bar R_{q+1}, \widetilde R_{q+1})\]
satisfies (\ref{pm-q}) with $q$ replaced by $q+1$.

The estimates (\ref{est-correct-1}) and (\ref{est-correct-2}) follow from Lemma \ref{le-correct}. Since $W_{I,q+1}$ is localized near frequency $\lambda_{q+1}$ in phase space, we can define
\begin{equation}\notag
\widetilde W_{I,q+1}=\nabla \Delta^{-1} \mathbb P_{\approx \lambda_{q+1}} (a_{I,q+1}e^{i\lambda_{q+1}\xi_I})
\end{equation}
and hence $W_{q+1}=\nabla\cdot \widetilde W_{q+1}$ with $\widetilde W_{q+1}=\sum_{I}\widetilde W_{I,q+1}$. Then the estimate (\ref{est-correct-3}) follows from Lemma \ref{le-amp1} with $k=r=0$. Regarding the advective derivative, we can rewrite 
\begin{equation}\notag
\begin{split}
&(\partial_t+T[M_{\e,q}]\cdot\nabla) \widetilde W_{I, q+1}\\
=& \left[ \frac{D_{M,q,\e}}{Dt}, \nabla\Delta^{-1}\mathbb P_{\approx \lambda_{q+1}}\right] (a_{I,q+1}e^{i\lambda_{q+1}\xi_I})+\nabla\Delta^{-1}\mathbb P_{\approx \lambda_{q+1}}\left(e^{i\lambda_{q+1}\xi_I}\frac{D_{M,q,\e}}{Dt}a_{I,q+1}\right).
\end{split}
\end{equation}
As a consequence, 
it follows from Lemma \ref{le-commu} with a suitable rescaling and Lemma \ref{le-amp1} with $k=0$ and $r=1$
\begin{equation}\notag
\begin{split}
\|(\partial_t+T[M_{\e,q}]\cdot\nabla) \widetilde W_{I, q+1}\|_{C^0}\lesssim&\ \lambda_q\delta_{q-1}^{\frac12}\lambda_{q+1}^{-1}\delta_{q}^{\frac12}+\lambda_{q+1}^{-1}\delta_{q}^{\frac12}\tau_{q}^{-1}\\
\lesssim&\ \lambda_{q+1}^{-1}\delta_{q}^{\frac12}\tau_{q}^{-1}
\end{split}
\end{equation}
provided 
\begin{equation}\label{tau-cond}
\tau_{q}\lesssim \lambda_q^{-1}\delta_{q-1}^{-\frac12}.
\end{equation}
Other terms in (\ref{est-correct-4}) can be estimated similarly.

By (\ref{induct-q-1}) and (\ref{est-correct-1}), we have
\begin{equation}\notag
\begin{split}
\|\nabla^k P_{q+1}\|_{C^0}\leq&\ \|\nabla^k P_{q}\|_{C^0} +\|\nabla^k W_{q+1}\|_{C^0}\\
\lesssim&\ \lambda_q^k\delta_{q-1}^{\frac12}+\lambda_{q+1}^k\delta_q^{\frac12}\\
\lesssim&\ \lambda_{q+1}^k\delta_q^{\frac12}
\end{split}
\end{equation}
since $b>1$ and $k\geq 1$. The estimate of $\nabla^k M_{q+1}, \nabla^k T[P_{q+1}], \nabla^k T[M_{q+1}]$ and hence (\ref{induct-q-1}) with $q$ replaced by $q+1$ follows analogously. 

Next we show (\ref{induct-q-2}) with $q$ replaced by $q+1$. We write
\begin{equation}\notag
\begin{split}
&(\partial_t+T[M_{q+1}]\cdot\nabla) T[M_{q+1}]\\
=&\ (\partial_t+T[M_{q}]\cdot\nabla) T[M_{q}]+T[W_{q+1}]\cdot\nabla T[M_q]+T[W_{q+1}]\cdot\nabla T[W_{q+1}]\\
&+(\partial_t+T[M_{q}]\cdot\nabla) T[W_{q+1}]
\end{split}
\end{equation}
and further decompose 
\begin{equation}\notag
\begin{split}
&(\partial_t+T[M_{q}]\cdot\nabla) T[W_{q+1}]\\
=&(\partial_t+T[M_{\e,q}]\cdot\nabla) T[W_{q+1}]+(T[M_{q}]-T[M_{\e,q}])\cdot\nabla T[W_{q+1}].
\end{split}
\end{equation}
Immediately it follows from the induction assumption (\ref{induct-q-2}) 
\begin{equation}\notag
\|(\partial_t+T[M_{q}]\cdot\nabla) T[M_{q}]\|_{C^0}\lesssim \lambda_q\delta_{q-1},
\end{equation}
and the estimate (\ref{est-correct-2}) 
\begin{equation}\notag
\|(\partial_t+T[M_{\e,q}]\cdot\nabla) T[W_{q+1}]\|_{C^0}\lesssim \tau_{q}^{-1}\delta_{q}^{\frac12}.
\end{equation}
The assumption (\ref{induct-q-1}) and estimate (\ref{est-correct-1}) together yield
\begin{equation}\notag
\|T[W_{q+1}]\cdot\nabla T[M_q]\|_{C^0}\leq \|T[W_{q+1}]\|_{C^0}\|\nabla T[M_q]\|_{C^0}
\lesssim \delta_q^{\frac12} \lambda_q \delta_{q-1}^{\frac12}.
\end{equation}
The estimate (\ref{est-correct-1}) also implies
\begin{equation}\notag
\|T[W_{q+1}]\cdot\nabla T[W_{q+1}]\|_{C^0}\leq \|T[W_{q+1}]\|_{C^0}\|\nabla T[W_{q+1}]\|_{C^0}
\lesssim \delta_q^{\frac12} \lambda_{q+1} \delta_{q}^{\frac12}.
\end{equation}
In view of Lemma \ref{le-molli} and (\ref{est-correct-1}) we have
\begin{equation}\notag
\begin{split}
\| (T[M_{q}]-T[M_{\e,q}])\cdot\nabla T[W_{q+1}] \|_{C^0}\leq& \|T[M_{q}]-T[M_{\e,q}]\|_{C^0}\|\nabla T[W_{q+1}]\|_{C^0}\\
\lesssim &\ \mu_q^{-L}\lambda_q^L\delta_{q-1}^{\frac12} \lambda_{q+1} \delta_{q}^{\frac12}.
\end{split}
\end{equation}
Summarizing the estimates above we obtain for $b>1$ and $0<\beta<1$
\begin{equation}\notag
\begin{split}
\|(\partial_t+T[M_{q+1}]\cdot\nabla) T[M_{q+1}]\|_{C^0}\lesssim&\ \lambda_q\delta_{q-1}+\tau_{q}^{-1}\delta_{q}^{\frac12}+ \delta_q^{\frac12} \lambda_q \delta_{q-1}^{\frac12}\\
&+\lambda_{q+1}\delta_{q}+\mu_q^{-L}\lambda_q^L\delta_{q-1}^{\frac12} \lambda_{q+1} \delta_{q}^{\frac12}\\
\lesssim&\ \lambda_{q+1}\delta_{q} 
\end{split}
\end{equation}
where we used $\mu_q=\lambda_{q+1}^{\frac{1}{L}}\lambda_q^{1-\frac{1}{L}}$ and supposed 
\begin{equation}\label{tau-cond2}
\tau_q^{-1}\leq \lambda_{q+1}\delta_q^{\frac12}.
\end{equation}
For $1\leq k\leq L-1$, higher order derivatives in (\ref{induct-q-2}) with $q$ replaced by $q+1$ can be estimated similarly. 

Next we establish the estimates for the new stress field.
Invoking $\widetilde \theta_q=\widetilde \theta_{q-1}$ and recalling (\ref{R-q1}), we have
\begin{equation}\notag
\begin{split}
\nabla\cdot \widetilde R_{q+1}=& \left(\partial_t+T[\widetilde \theta_{\e, q-1}]\cdot\nabla \right)W_{q+1}+\nu\Lambda^\gamma W_{q+1}
+T[W_{q+1}]\cdot\nabla\widetilde\theta_{\e, q-1}\\
&+\nabla\cdot\left(c_{1,q}A_1+R^*_q-2T[W_{q+1}] W_{q+1} \right)\\
=&\left(\partial_t+T[\widetilde \theta_{\e, q-1}]\cdot\nabla \right)W_{q+1}+\nu\Lambda^\gamma W_{q+1}
+T[W_{q+1}]\cdot\nabla\widetilde\theta_{\e, q-1}\\
&+\nabla\cdot\left(c_{\e,1,q}A_1+R^*_{\e, q}-2T[W_{q+1}] W_{q+1} \right) \\
&+\nabla\cdot \left( (T[\widetilde \theta_{q-1}]-T[\widetilde \theta_{\e, q-1}])W_{q+1}+T[W_{q+1}](\widetilde \theta_{q-1}-\widetilde \theta_{\e, q-1}) \right.\\
&\left. +(c_{1,q}-c_{\e,1,q})A_1+(R^*_q-R^*_{\e, q})\right)\\
=&: \nabla\cdot R_T+\nabla\cdot R_D+\nabla\cdot R_N+\nabla\cdot R_O+\nabla\cdot R_M
\end{split}
\end{equation}
with $\widetilde \theta_{\e,q-1}= P_{\e,q-1}- M_{\e,q-1}$.

\medskip

\textbf{Estimates of $R_T$:} 
Note $W_{q+1}$ is localized to frequency $\approx \lambda_{q+1}$ in Fourier space. 
Therefore
we can find $R_T$ such that
\[R_T=\nabla\Delta^{-1} \mathbb P_{\approx \lambda_{q+1}} \left[(\partial_t+T[\widetilde \theta_{\e, q-1}]\cdot\nabla) W_{q+1}\right].\]
As a consequence we obtain
\begin{equation}\notag
\|R_T\|_{C^0}\lesssim \lambda_{q+1}^{-1} \|(\partial_t+T[\widetilde \theta_{\e, q-1}]\cdot\nabla)W_{q+1}\|_{C^0}.
\end{equation}
Since \[T[\widetilde \theta_{\e, q-1}]=T[P_{\e, q-1}]-T[M_{\e, q-1}],\]
applying Lemma \ref{le-correct} with $q'=q$, $r=1$ and $k=0$ we have
\begin{equation}\notag
\begin{split}
&\|(\partial_t+T[\widetilde \theta_{\e, q-1}]\cdot\nabla)W_{q+1}\|_{C^0}\\
\leq& \|\left(\partial_t+T[P_{\e, q-1}]\cdot\nabla \right)W_{q+1}\|_{C^0}+\|\left(\partial_t+T[M_{\e, q-1}]\cdot\nabla \right)W_{q+1}\|_{C^0}\\
\lesssim&\ \tau_{q}^{-1}\delta_q^{\frac12}.
\end{split}
\end{equation}
Therefore, we conclude
\begin{equation}\label{est-RT}
\|R_T\|_{C^0}\lesssim \lambda_{q+1}^{-1} \tau_{q}^{-1}\delta_q^{\frac12}.
\end{equation}

\medskip

\textbf{Estimates of $R_D$:} It is obvious that there exists $R_D$ satisfying
\[R_D=\nu\nabla\Delta^{-1}\Lambda^\gamma W_{q+1}.\]
It follows from (\ref{est-correct-1}) that
\begin{equation}\label{est-RD}
\|R_D\|_{C^0}\lesssim \lambda_{q+1}^{-1+\gamma} \|W_{q+1}\|_{C^0}\lesssim \lambda_{q+1}^{-1+\gamma}\delta_q^{\frac12}.
\end{equation}

\medskip

\textbf{Estimates of $R_N$:} Again, due to the frequency localization property of $\widetilde \theta_{\e, q-1}$ and $W_{q+1}$, we can define
\begin{equation}\notag
R_N=\nabla\Delta^{-1}\mathbb P_{\approx \lambda_{q+1}}[ T[W_{q+1}]\cdot\nabla\widetilde\theta_{\e, q-1}].
\end{equation}
Applying (\ref{induct-q-7}) and (\ref{est-correct-1}) gives
\begin{equation}\label{est-RN}
\begin{split}
\|R_N\|_{C^0}\lesssim &\ \lambda_{q+1}^{-1}\|T[W_{q+1}]\|_{C^0} \|\nabla\widetilde\theta_{\e, q-1}\|_{C^0} \\
\lesssim &\ \lambda_{q+1}^{-1}\delta_q^{\frac12}\lambda_{q-1}\delta_{q-2}^{\frac12}.
\end{split}
\end{equation}

\medskip

\textbf{Estimates of $R_O$:} 
It follows from Microlocal Lemma \ref{le-mic} that
\begin{equation}\notag
\begin{split}
W_{I,q+1}=&\ \mathbb P_{\approx \lambda_{q+1}}[a_{I,q+1}e^{i\lambda_{q+1}\xi_I}]\\
=&\ a_{I,q+1}\eta_{\lambda_{q+1}}(\lambda_{q+1}\nabla\xi_{I})e^{i\lambda_{q+1}\xi_I}+\delta a_{I,q+1}e^{i\lambda_{q+1}\xi_I}\\
=&\ a_{I,q+1}e^{i\lambda_{q+1}\xi_I}+\delta a_{I,q+1}e^{i\lambda_{q+1}\xi_I}
\end{split}
\end{equation}
and
\begin{equation}\notag
\begin{split}
T[W_{I,q+1}]=&\ T \mathbb P_{\approx \lambda_{q+1}}[a_{I,q+1}e^{i\lambda_{q+1}\xi_I}]\\
=&\ a_{I,q+1}m(\lambda_{q+1}\nabla\xi_{I})\eta_{\lambda_{q+1}}(\lambda_{q+1}\nabla\xi_{I})e^{i\lambda_{q+1}\xi_I}+\delta u_{I,q+1}e^{i\lambda_{q+1}\xi_I}\\
=&\ a_{I,q+1}m(\nabla\xi_{I})e^{i\lambda_{q+1}\xi_I}+\delta u_{I,q+1}e^{i\lambda_{q+1}\xi_I}.
\end{split}
\end{equation}
Consequently we compute
\begin{equation}\notag
\begin{split}
T[W_{q+1}] \cdot\nabla W_{q+1}=&\ \frac12\cdot\nabla\sum_{I\in\Omega} \left(T[W_{I,q+1}] W_{\bar I, q+1}+T[W_{\bar I,q+1}] W_{I, q+1}\right)\\
&+\sum_{J\neq \bar I} T[W_{J, q+1}] \cdot\nabla W_{I, q+1}\\
\end{split}
\end{equation}
with 
\begin{equation}\notag
\begin{split}
&\frac12\sum_{I\in\Omega} \left(T[W_{I,q+1}] W_{\bar I, q+1}+T[W_{\bar I,q+1}] W_{I, q+1}\right)\\
=&\ \frac12\sum_{I\in\Omega} a_{I,q+1}^2\left( m(\nabla\widehat\xi_{I})+m(-\nabla\widehat\xi_{I})\right)\\
&+\frac12\sum_{I\in\Omega} a_{I,q+1}^2\left( m(\nabla\xi_{I})-m(\nabla\widehat\xi_{I})+m(-\nabla\xi_{I})-m(-\nabla\widehat\xi_{I})\right)\\
&+\frac12\sum_{I\in\Omega} \left( a_{I,q+1}\delta u_{\bar I, q+1}+a_{I,q+1}m(\nabla\xi_I)\delta a_{\bar I,q+1}+a_{\bar I,q+1}\delta u_{I, q+1}\right.\\
&\left. +a_{\bar I, q+1}m(\nabla\xi_{\bar I}) \delta a_{I,q+1}-\delta a_{\bar I, q+1}\delta u_{I, q+1}-\delta a_{I, q+1}\delta u_{\bar I, q+1}\right)\\
=&: \frac12\sum_{I\in\Omega} a_{I,q+1}^2\left( m(\nabla\widehat\xi_{I})+m(-\nabla\widehat\xi_{I})\right)+R_{O,1}+R_{O,2}.
\end{split}
\end{equation}
Since 
\begin{equation}\notag
\nabla W_{I,q+1}= i\lambda_{q+1}\nabla\xi_I\mathbb P_{\approx \lambda_{q+1}}[a_{I,q+1}e^{i\lambda_{q+1}\xi_I}]
+\mathbb P_{\approx \lambda_{q+1}}[\nabla a_{I,q+1}e^{i\lambda_{q+1}\xi_I}],
\end{equation}
we further compute
\begin{equation}\notag
\begin{split}
&\sum_{J\neq \bar I} T[W_{J, q+1}]\cdot\nabla W_{I, q+1}\\
=&\sum_{J\neq \bar I} i\lambda_{q+1}T[W_{J, q+1}]\cdot \nabla\xi_I\mathbb P_{\approx \lambda_{q+1}}[a_{I,q+1}e^{i\lambda_{q+1}\xi_I}]\\
&+\sum_{J\neq \bar I}T[W_{J, q+1}]\cdot \mathbb P_{\approx \lambda_{q+1}}[\nabla a_{I,q+1}e^{i\lambda_{q+1}\xi_I}]\\
=&\sum_{J\neq \bar I} i\lambda_{q+1}a_{J,q+1}m(\nabla\xi_{J})e^{i\lambda_{q+1}\xi_J}\cdot \nabla\xi_I\mathbb P_{\approx \lambda_{q+1}}[a_{I,q+1}e^{i\lambda_{q+1}\xi_I}]\\
&+\sum_{J\neq \bar I} i\lambda_{q+1}\delta u_{J,q+1}e^{i\lambda_{q+1}\xi_J}\cdot \nabla\xi_I\mathbb P_{\approx \lambda_{q+1}}[a_{I,q+1}e^{i\lambda_{q+1}\xi_I}]\\
&+\sum_{J\neq \bar I}T[W_{J, q+1}]\cdot \mathbb P_{\approx \lambda_{q+1}}[\nabla a_{I,q+1}e^{i\lambda_{q+1}\xi_I}].
\end{split}
\end{equation}
Note that
\begin{equation}\notag
\begin{split}
m(\nabla\xi_J)\cdot\nabla \xi_I
=&\ m(\nabla\xi_{J,in})\cdot\nabla \xi_{I,in}+m(\nabla\xi_{J,in})\cdot(\nabla \xi_{I}-\nabla \xi_{I,in})\\
&+(m(\nabla\xi_J)-m(\nabla\xi_{J,in}))\cdot\nabla \xi_I\\
=&\ m(\nabla\xi_{J,in})\cdot(\nabla \xi_{I}-\nabla \xi_{I,in})\\
&+(m(\nabla\xi_J)-m(\nabla\xi_{J,in}))\cdot\nabla \xi_I\\
\end{split}
\end{equation}
since 
\[m(\nabla\xi_{J,in})\cdot\nabla \xi_{I,in}=m(\pm\nabla\xi_{I,in})\cdot\nabla \xi_{I,in}=0.\]
Therefore we have
\begin{equation}\notag
\begin{split}
&\sum_{J\neq \bar I} T[W_{J, q+1}]\cdot\nabla W_{I, q+1}\\
=&\sum_{J\neq \bar I} i\lambda_{q+1}a_{J,q+1}m(\nabla\xi_{J,in})\cdot(\nabla \xi_{I}-\nabla \xi_{I,in})e^{i\lambda_{q+1}\xi_J}\mathbb P_{\approx \lambda_{q+1}}[a_{I,q+1}e^{i\lambda_{q+1}\xi_I}]\\
&+\sum_{J\neq \bar I} i\lambda_{q+1}a_{J,q+1}(m(\nabla\xi_J)-m(\nabla\xi_{J,in}))\cdot\nabla \xi_I e^{i\lambda_{q+1}\xi_J}\mathbb P_{\approx \lambda_{q+1}}[a_{I,q+1}e^{i\lambda_{q+1}\xi_I}]\\
&+\sum_{J\neq \bar I} i\lambda_{q+1}\delta u_{J,q+1}e^{i\lambda_{q+1}\xi_J}\cdot \nabla\xi_I\mathbb P_{\approx \lambda_{q+1}}[a_{I,q+1}e^{i\lambda_{q+1}\xi_I}]\\
&+\sum_{J\neq \bar I}T[W_{J, q+1}]\cdot \mathbb P_{\approx \lambda_{q+1}}[\nabla a_{I,q+1}e^{i\lambda_{q+1}\xi_I}]\\
=&: \nabla\cdot R_{O,3}+\nabla\cdot R_{O,4}+\nabla\cdot R_{O,5}+\nabla\cdot R_{O,6}.
\end{split}
\end{equation}
Summarizing the analysis above we obtain
\begin{equation}\notag
\begin{split}
T[W_{q+1}]\cdot\nabla W_{q+1}=&\ \frac12\nabla\cdot \sum_{I\in\Omega} a_{I,q+1}^2\left( m(\nabla\xi_{I,in})+m(-\nabla\xi_{I,in})\right)+\nabla\cdot R_{O,1}\\
&+\nabla\cdot R_{O,2}+\nabla\cdot R_{O,3}+\nabla\cdot R_{O,4}+\nabla\cdot R_{O,5}+\nabla\cdot R_{O,6}.
\end{split}
\end{equation}
According to the choice of $a_{I, q+1}$ in (\ref{amplitude}), we have that
\begin{equation}\notag
\begin{split}
\nabla\cdot R_O=&\ \nabla\cdot\left(c_{\e, 1,q}A_1+R^*_{\e,q}-2T[W_{q+1}] W_{q+1} \right)\\
=&\ \nabla\cdot \left( c_{2, q+1} A_2-2R_{O,1} -2R_{O,2}-2R_{O,3}-2R_{O,4}-2R_{O,5}-2R_{O,6}\right).
\end{split}
\end{equation}
Now we estimate the error terms above. It follows from the definition of $c_{2, q+1} $ and the assumption (\ref{induct-q-5}) that
\begin{equation}\notag
\|c_{2, q+1}\|_{C^0}\lesssim \|R^*_{\e, q}\|_{C^0}\lesssim \|R^*_{q}\|_{C^0}\lesssim \delta_{q+1}.
\end{equation}
Applying Lemma \ref{le-phase} and \ref{le-amp1}, noticing that $\xi_I$ is advected by $\widetilde u_q=\widetilde u_{q-1}$, leads to
\begin{equation}\notag
\begin{split}
\|R_{O,1}\|_{C^0}\lesssim & \sum_{I\in\Omega} \|a_{I,q+1}\|_{C^0}^2|\nabla\xi_I-\nabla\widehat\xi_{I}|\\
\lesssim & \sum_{I\in\Omega} \|a_{I,q+1}\|_{C^0}^2\lambda_{q-1}\tau_{q}\|\widetilde u_{q-1}\|_{C^0}\\
\lesssim&\ \lambda_{q-1}\tau_{q}\delta_{q-2}^{\frac12}\delta_q.
\end{split}
\end{equation}
Applying Lemma \ref{le-amp1} and Lemma \ref{le-amp2}
\begin{equation}\notag
\begin{split}
\|R_{O,2}\|_{C^0}\lesssim & \sum_{I\in\Omega} \|a_{I,q+1}\|_{C^0} \|\delta u_{\bar I,q+1}\|_{C^0}+  \sum_{I\in\Omega}\|a_{I,q+1}\|_{C^0} \|\delta a_{\bar I,q+1}\|_{C^0} \\
&+ \sum_{I\in\Omega}\|\delta a_{I,q+1}\|_{C^0} \|\delta u_{\bar I,q+1}\|_{C^0}  \\
\lesssim&\ \delta_{q}^{\frac12}\lambda_{q+1}^{-1}\lambda_q\delta_q^{\frac12}.
\end{split}
\end{equation}
We observe that due to the frequency support of $W_{q+1}$ and $T[W_{q+1}]$, we can define 
\begin{equation}\notag
R_{O,3}+R_{O,4}+R_{O,5}+R_{O,6}=\sum_{J\neq \bar I} \nabla\Delta^{-1}\mathbb P_{\approx \lambda_{q+1}}\left[T[W_{J, q+1}]\cdot\nabla W_{I, q+1}\right].
\end{equation}
Therefore we deduce from Lemma \ref{le-phase} and Lemma \ref{le-amp1}
\begin{equation}\notag
\begin{split}
\|R_{O,3}\|_{C^0}\lesssim & \sum_{J\neq \bar I}\|a_{J,q+1}\|_{C^0}\|a_{I,q+1}\|_{C^0}\|m(\nabla \xi_{I,in})\|_{C^0}|\nabla\xi_{I}-\nabla\xi_{I,in}|\\
\lesssim & \sum_{J\neq \bar I}\|a_{J,q+1}\|_{C^0}\|a_{I,q+1}\|_{C^0}\lambda_{q-1}\tau_{q}\|\widetilde u_{q-1}\|_{C^0}\\
\lesssim&\ \lambda_{q-1}\tau_{q} \delta_{q-2}^{\frac12}\delta_q
\end{split}
\end{equation}
\begin{equation}\notag
\begin{split}
\|R_{O,4}\|_{C^0}\lesssim & \sum_{J\neq \bar I}\|a_{J,q+1}\|_{C^0}\|a_{I,q+1}\|_{C^0}\|m(\nabla\xi_{J})-m(\nabla\xi_{J,in})\|_{C^0}\|\nabla\xi_I\|_{C^0}\\
\lesssim & \sum_{J\neq \bar I}\|a_{J,q+1}\|_{C^0}\|a_{I,q+1}\|_{C^0}|\nabla\xi_{I}-\nabla\xi_{I,in}|\\
\lesssim&\ \lambda_{q-1}\tau_{q} \delta_{q-2}^{\frac12} \delta_q
\end{split}
\end{equation}
where we used the fact 
$ \|m(\nabla\xi_{J})-m(\nabla\xi_{J,in})\|_{C^0}\lesssim \|\nabla\xi_{J}-\nabla\xi_{J,in}\|_{C^0}$;
\begin{equation}\notag
\begin{split}
\|R_{O,5}\|_{C^0}\lesssim & \sum_{J\neq \bar I}\|\delta u_{J,q+1}\|_{C^0}\|a_{I,q+1}\|_{C^0}\|\nabla\xi_I\|_{C^0}\\
\lesssim&\ \lambda_{q+1}^{-1}\lambda_q \delta_{q}^{\frac12}\delta_q^{\frac12}
\end{split}
\end{equation}
using (\ref{est-correct-1}) and Lemma \ref{le-amp2}, and 
\begin{equation}\notag
\begin{split}
\|R_{O,6}\|_{C^0}\lesssim & \sum_{J\neq \bar I}\lambda_{q+1}^{-1}\|T[W_{J,q+1}]\|_{C^0}\|\nabla a_{I,q+1}\|_{C^0}\\
\lesssim&\ \lambda_{q+1}^{-1}\delta_{q}^{\frac12}\lambda_q \delta_{q}^{\frac12}.
\end{split}
\end{equation}
Summarizing the estimates above gives
\begin{equation}\label{est-RO1}
\|R_{O,1}\|_{C^0}+\|R_{O,3}\|_{C^0}+\|R_{O,4}\|_{C^0}\lesssim  \lambda_{q-1}\tau_{q} \delta_{q-2}^{\frac12} \delta_q,
\end{equation}
\begin{equation}\label{est-RO2}
\|R_{O,2}\|_{C^0}+\|R_{O,5}\|_{C^0}+\|R_{O,6}\|_{C^0}\lesssim  \lambda_{q+1}^{-1}\lambda_q\delta_{q}.
\end{equation}

\medskip

\textbf{Estimates of $R_M$:} 
It follows from the fact $\widetilde \theta_{q}=\widetilde \theta_{q-1}$, Lemma \ref{le-molli} and estimate (\ref{est-correct-1}) that
\begin{equation}\notag
\begin{split}
&\|(T[\widetilde \theta_{q-1}]-T[\widetilde \theta_{\e, q-1}])W_{q+1}\|_{C^0}+\|T[W_{q+1}](\widetilde \theta_{q-1}-\widetilde \theta_{\e, q-1})\|_{C^0}\\
=&\ \|(T[\widetilde \theta_{q}]-T[\widetilde \theta_{\e, q}])W_{q+1}\|_{C^0}+\|T[W_{q+1}](\widetilde \theta_{q}-\widetilde \theta_{\e, q})\|_{C^0}\\
\lesssim&\ \delta_q^{\frac12}\mu_{q}^{-L}\lambda_{q}^{L}\delta_{q-1}^{\frac12}\\
\lesssim&\ \delta_{q}^{\frac12}\delta_{q-1}^{\frac12}\lambda_{q+1}^{-1}\lambda_q.
\end{split}
\end{equation}
By Lemma \ref{le-molli-est}, we have
\begin{equation}\notag
\|(c_{1, q, \e}-c_{1,q})A_1\|_{C^0}+\|R^*_{q,\e}-R^*_{q}\|_{C^0}\lesssim \delta_{q}^{\frac12}\delta_{q-1}^{\frac12}\lambda_{q+1}^{-1}\lambda_q.
\end{equation}
Therefore 
\begin{equation}\label{est-RM}
\|R_M\|_{C^0}\lesssim \delta_{q}^{\frac12}\delta_{q-1}^{\frac12}\lambda_{q+1}^{-1}\lambda_q.
\end{equation}

In summary, the new stress error can be written as
\begin{equation}\notag
\widetilde R_{q+1}=c_{2,q+1}A_2+R^*_{q+1}
\end{equation}
with 
\begin{equation}\notag
\|c_{2,q+1}\|_{C^0}\lesssim \delta_{q+1} 
\end{equation}
and 
\[R^*_{q+1}=R_T+R_N+R_D+R_M-2(R_{O,1}+...+R_{O,6}).\]
Thus (\ref{induct-q-3}) is satisfied with $c_{1,q}$ replaced by $c_{2,q+1}$. To show (\ref{induct-q-5}) with $q$ replaced by $q+1$, we just need to show $\|R^*_{q+1}\|_{C^0}\leq \delta_{q+2}$. We choose $\tau_q=\lambda_{q+1}^{-\frac12}\lambda_{q-1}^{-\frac12}\delta_q^{-\frac14}\delta_{q-2}^{-\frac14}$ to optimize the two estimates (\ref{est-RT}) and (\ref{est-RO1})  such that
\begin{equation}\notag
\|R_T\|_{C^0}+\|R_{O,1}\|_{C^0}+\|R_{O,3}\|_{C^0}+\|R_{O,4}\|_{C^0}\lesssim \lambda_{q+1}^{-\frac12}\lambda_{q-1}^{\frac12}\delta_q^{\frac34}\delta_{q-2}^{\frac14}.
\end{equation}
In view of the last estimate, (\ref{est-RN}), (\ref{est-RD}), (\ref{est-RO2}) and (\ref{est-RM}), we impose
\begin{equation}\notag
\begin{cases}
\lambda_{q+1}^{-\frac12}\lambda_{q-1}^{\frac12}\delta_q^{\frac34}\delta_{q-2}^{\frac14}\lesssim \delta_{q+2}\\
\lambda_{q+1}^{-1}\lambda_{q-1}\delta_q^{\frac12}\delta_{q-2}^{\frac12}\lesssim \delta_{q+2}\\
\lambda_{q+1}^{-1+\gamma}\delta_q^{\frac12}\lesssim \delta_{q+2}\\
\lambda_{q+1}^{-1}\lambda_q\delta_q\lesssim \delta_{q+2}\\
\lambda_{q+1}^{-1}\lambda_q\delta_q^{\frac12}\delta_{q-1}^{\frac12}\lesssim \delta_{q+2}\\
\end{cases}
\end{equation}
which are satisfied provided $b=1^+$ (close enough to 1 from the right), $\beta<\frac25$ and for $\gamma\geq 0$ satisfying $\beta<\frac{2b(1-\gamma)}{2b^2-1}$. Recall $\theta$ is in $C_t^0C_x^\alpha$ with $\alpha<\frac{\beta}{2b}<\frac15$. The conditions $\beta<\frac{2b(1-\gamma)}{2b^2-1}$ and $\alpha<\frac{\beta}{2b}$ together imply $0\leq \gamma<1-\alpha$.

\medskip

\textbf{Estimates of higher order spatial derivatives:} For $\nabla^k c_{2,q+1}$, each derivative cost is $\lesssim \lambda_{q+1}$. Combining with the $C^0$ estimate of $c_{2,q+1}$, we have 
\begin{equation}\notag
\|\nabla^k c_{2,q+1}\|_{C^0}\lesssim \lambda_{q+1}^k\delta_{q+1}.
\end{equation}
We observe that $R_T, R_D, R_N$, $R_{O,3}$, $R_{O,4}$,$R_{O,5}$ and $R_{O,6}$ are localized in Fourier space near frequency $\lambda_{q+1}$. Hence each spatial derivative of them costs at most $\approx \lambda_{q+1}$. From Lemma \ref{le-molli}, we know 
\[ \|\widetilde \theta_{q-1}-\widetilde \theta_{q-1, \e} \|_{C^0} \lesssim \lambda_q^{-1}\lambda_{q-1} \delta_{q-2}^{\frac12}.\]
We also have the estimate from (\ref{induct-q-7})
\[ \|\nabla \widetilde \theta_{q-1}\|_{C^0}+\|\nabla \widetilde \theta_{q-1,\e}\|_{C^0}\lesssim \lambda_{q-1}\delta_{q-2}^{\frac12}. \]
Therefore one spatial derivative cost of $\widetilde \theta_{q-1}-\widetilde \theta_{q-1, \e}$ and $T[\widetilde \theta_{q-1}]-T[\widetilde \theta_{q-1, \e}]$ is at most $\approx \lambda_q$. On the other hand, one spatial derivative cost of $W_{q+1}$ and $T[W_{q+1}]$ is at most $\approx \lambda_{q+1}$. Hence one spatial derivative cost of $R_M$ is at most $\approx \lambda_{q+1}$. Regarding $R_{O,1}$, we compare the estimates from Lemma \ref{le-phase}
\begin{equation}\notag
\begin{split}
\|\nabla^2\xi_I(\Phi_q(x,s,t))\|_{C^0}\lesssim&\ \lambda_q,\\
\|\nabla\xi_I(\Phi_q(x,s,t))-\nabla\xi_{I,in}\|_{C^0}\lesssim&\ \tau_{q}\lambda_q\delta_{q-1}^{\frac12}.
\end{split}
\end{equation}
Thus one spatial derivative of $\nabla\xi_I(\Phi_q(x,s,t))-\nabla\xi_{I,in}$ is at most $\approx \tau_{q}^{-1}\delta_{q-1}^{-\frac12}\lesssim \lambda_{q+1}$.
The analysis shows that 
\begin{equation}\notag
\|\nabla^k R^*_{q+1}\|_{C^0}\lesssim \lambda_{q+1}^k\delta_{q+2}, \ \ 0\leq k\leq L.
\end{equation}

\medskip

\textbf{Estimates of advective derivative:} 
First we rewrite
\begin{equation}\notag
\partial_t+T[M_{q+1}]\cdot\nabla=\partial_t+T[M_{q, \e}]\cdot\nabla+(T[M_{q}]-T[M_{q, \e}])\cdot\nabla+T[W_{q+1}]\cdot\nabla
\end{equation}
and 
\begin{equation}\notag
\partial_t+T[\theta_{q+1}]\cdot\nabla=\partial_t+T[\theta_{q, \e}]\cdot\nabla+(T[\theta_{q}]-T[\theta_{q, \e}])\cdot\nabla+T[W_{q+1}]\cdot\nabla.
\end{equation}
For $R_D$: 
\begin{equation}\notag
\begin{split}
&\|\nabla^k(\partial_t+T[M_{q+1}]\cdot\nabla)R_D\|_{C^0}\\
\lesssim&\ \|\nabla^k(\partial_t+T[M_{q,\e}]\cdot\nabla)R_{D}\|_{C^0}+\|\nabla^k\left((T[M_{q}]-T[M_{q,\e}])\cdot\nabla\right) R_{D}\|_{C^0}\\
&+\|\nabla^k(T[W_{q+1}]\cdot\nabla) R_{D}\|_{C^0}\\
\lesssim&\ \lambda_{q+1}^{k+1}\delta_{q}^{\frac12}\delta_{q+1}.
\end{split}
\end{equation}
In view of the definition of $c_{2,q+1}$ in (\ref{c2-q1}), we observe
\begin{equation}\notag
\begin{split}
&\|\nabla^k(\partial_t+T[M_{q+1}]\cdot\nabla)c_{2,q+1}\|_{C^0}\\
\lesssim&\ \|\nabla^k(\partial_t+T[M_{q+1}]\cdot\nabla)R^*_{q,\e}\|_{C^0}\\
\lesssim&\ \|\nabla^k(\partial_t+T[M_{q+1}]\cdot\nabla)R^*_{q}\|_{C^0}\\
\lesssim&\ \|\nabla^k(\partial_t+T[M_{q,\e}]\cdot\nabla)R^*_{q}\|_{C^0}+\|\nabla^k\left((T[M_{q}]-T[M_{q,\e}])\cdot\nabla\right) R^*_{q}\|_{C^0}\\
&+\|\nabla^k(T[W_{q+1}]\cdot\nabla) R^*_{q}\|_{C^0}\\
\lesssim&\ \lambda_q^{k+1} \delta_{q-1}^{\frac12}\delta_{q+1}+\lambda_{q+1}^k\delta_{q}^{\frac12}\lambda_q\delta_{q+1}\\
\lesssim&\ \lambda_{q+1}^{k+1}\delta_{q}^{\frac12}\delta_{q+1}.
\end{split}
\end{equation}

\begin{Lemma}\label{le-phase-high}
Let $k\geq 0$ and $0\leq r\leq 2$, we have
\begin{equation}\notag
\|\nabla^k\left(\frac{D_q}{Dt} \right)^r(\nabla\xi_I-\nabla\xi_{I,in})\|_{C^0}\lesssim \lambda_{q+1}^k\tau_{q}^{-r}  \tau_{q+1}\lambda_q\delta_{q-1}^{\frac12}.
\end{equation}
\end{Lemma}

\begin{Lemma}\label{le-m-diff}
Let $k\geq 0$ and $0\leq r\leq 2$, we also have
\begin{equation}\notag
\|\nabla^k\left(\frac{D_q}{Dt} \right)^r(m(\nabla\xi_I)-m(\nabla\xi_{I,in}))\|_{C^0}\lesssim \lambda_{q+1}^k\tau_{q}^{-r}  \tau_{q}\lambda_q\delta_{q-1}^{\frac12}.
\end{equation}
\end{Lemma}
\pf
Note
\begin{equation}\notag
m(\nabla\xi_I)-m(\nabla\xi_{I,in})=\left(\nabla\xi_I-\nabla\xi_{I,in}\right)\int_0^1 \partial_a m\left((1-s)\nabla\xi_{I,in}+s\nabla\xi_I\right)\, ds.
\end{equation}
The estimate follows from Lemma \ref{le-phase-high}.
\cbdu

Combining Lemma \ref{le-amp1}, Lemma \ref{le-amp2}, Lemma \ref{le-phase-high}, Lemma \ref{le-m-diff}, we have
\begin{equation}\notag
\begin{split}
&\|\nabla^k(\partial_t+T[M_{q+1}]\cdot\nabla)R_{O,1}\|_{C^0}+\|\nabla^k(\partial_t+T[M_{q+1}]\cdot\nabla)R_{O,2}\|_{C^0}\\
\lesssim&\ \|\nabla^k(\partial_t+T[M_{q,\e}]\cdot\nabla)R_{O,1}\|_{C^0}+\|\nabla^k\left((T[M_{q}]-T[M_{q,\e}])\cdot\nabla\right) R_{O,1}\|_{C^0}\\
&+\|\nabla^k(T[W_{q+1}]\cdot\nabla) R_{O,1}\|_{C^0}+\|\nabla^k(\partial_t+T[M_{q,\e}]\cdot\nabla)R_{O,2}\|_{C^0}\\
&+\|\nabla^k\left((T[M_{q}]-T[M_{q,\e}])\cdot\nabla\right) R_{O,2}\|_{C^0}
+\|\nabla^k(T[W_{q+1}]\cdot\nabla) R_{O,2}\|_{C^0}\\
\lesssim&\ \lambda_{q+1}^{k+1}\delta_{q}^{\frac12}\delta_{q+2}.
\end{split}
\end{equation}

For $R_T$, recall 
\[R_T=\nabla\Delta^{-1}\mathbb P_{\approx \lambda_{q+1}} [(\partial_t+T[\widetilde \theta_{q-1}]\cdot\nabla)W_{q+1}].\]
Then
\begin{equation}\notag
\begin{split}
&\|(\partial_t+T[\theta_{q+1}]\cdot\nabla) R_T\|_{C^0}\\
\lesssim& \ \|(\partial_t+T[\theta_{q,\e}]\cdot\nabla) R_T\|_{C^0}+\|(T[\theta_q]-T[\theta_{q,\e}])\cdot\nabla R_T\|_{C^0}\\
&+\|T[W_{q+1}]\cdot\nabla R_T\|_{C^0}.
\end{split}
\end{equation}
Denote 
\[\frac{D_{q,\theta,\e}}{Dt}= \partial_t+T[\theta_{q,\e}]\cdot\nabla, \ \ \frac{\widetilde D_{q,\theta}}{Dt}= \partial_t+T[\widetilde \theta_{q}]\cdot\nabla\]
For the first term, we apply the commutator 
\begin{equation}\notag
\begin{split}
(\partial_t+T[\theta_{q,\e}]\cdot\nabla) R_T=&\left[\frac{D_{q,\theta,\e}}{Dt}, \nabla\Delta^{-1}\mathbb P_{\approx \lambda_{q+1}}\right] \frac{\widetilde D_{q-1,\theta}}{Dt} W_{q+1}\\
&+\nabla\Delta^{-1}\mathbb P_{\approx \lambda_{q+1}}\frac{D_{q,\theta,\e}}{Dt} \frac{\widetilde D_{q-1,\theta}}{Dt} W_{q+1}.
\end{split}
\end{equation}
Hence applying Lemma \ref{le-commu}, Lemma \ref{le-correct}
\begin{equation}\notag
\begin{split}
\|(\partial_t+T[\theta_{q,\e}]\cdot\nabla) R_T\|_{C^0}\leq&\ \|\left[\frac{D_{q,\theta,\e}}{Dt}, \nabla\Delta^{-1}\mathbb P_{\approx \lambda_{q+1}}\right] \frac{\widetilde D_{q-1,\theta}}{Dt} W_{q+1}\|_{C^0}\\
&+\|\nabla\Delta^{-1}\mathbb P_{\approx \lambda_{q+1}}\frac{D_{q,\theta,\e}}{Dt} \frac{\widetilde D_{q-1,\theta}}{Dt} W_{q+1}\|_{C^0}\\
\lesssim&\ \lambda_q\delta_{q-1}^{\frac12} \tau_{q-1}^{-1} \delta_q^{\frac12}+\lambda_{q+1}^{-1}\tau_{q}^{-1} \tau_{q-1}^{-1}\delta_q^{\frac12}.
\end{split}
\end{equation}
Applying Lemma \ref{le-molli}, estimate (\ref{est-correct-1}) and the spatial derivative established earlier
\begin{equation}\notag
\begin{split}
\|(T[\theta_q]-T[\theta_{q,\e}])\cdot\nabla R_T\|_{C^0}\lesssim&\ \|T[\theta_q]-T[\theta_{q,\e}]\|_{C^0}\|\nabla R_T\|_{C^0}\\
\lesssim&\ \lambda_{q+1}^{-1}\lambda_q\delta_{q-1}^{\frac12} \lambda_{q+1}\delta_{q+2},\\
\end{split}
\end{equation}
\begin{equation}\notag
\begin{split}
\|T[W_{q+1}]\cdot\nabla R_T\|_{C^0}\lesssim&\ \|T[W_{q+1}]\|_{C^0}\|\nabla R_T\|_{C^0}\\
\lesssim&\ \delta_{q}^{\frac12} \lambda_{q+1}\delta_{q+2}.
\end{split}
\end{equation}
Summarizing we have
\begin{equation}\notag
\begin{split}
\|(\partial_t+T[\theta_{q+1}]\cdot\nabla) R_T\|_{C^0}\lesssim&\ \lambda_q\delta_{q-1}^{\frac12} \tau_{q-1}^{-1} \delta_q^{\frac12}+\lambda_{q+1}^{-1}\tau_{q}^{-1} \tau_{q-1}^{-1}\delta_q^{\frac12}\\
&+\lambda_{q+1}^{-1}\lambda_q\delta_{q-1}^{\frac12} \lambda_{q+1}\delta_{q+2}+ \delta_{q}^{\frac12} \lambda_{q+1}\delta_{q+2}.
\end{split}
\end{equation}
The advective derivative of $R_N$ can be handled similarly.
Recall 
\[R_N=\nabla\Delta^{-1}\mathbb P_{\approx \lambda_{q+1}} [T[W_{q+1}]\cdot\nabla \widetilde \theta_{q-1}].\]
\begin{equation}\notag
\begin{split}
&\|(\partial_t+T[\theta_{q+1}]\cdot\nabla) R_N\|_{C^0}\\
\lesssim& \ \|(\partial_t+T[\theta_{q,\e}]\cdot\nabla) R_N\|_{C^0}+\|(T[\theta_q]-T[\theta_{q,\e}])\cdot\nabla R_N\|_{C^0}\\
&+\|T[W_{q+1}]\cdot\nabla R_N\|_{C^0}.
\end{split}
\end{equation}
\begin{equation}\notag
\begin{split}
(\partial_t+T[\theta_{q,\e}]\cdot\nabla) R_N=&\left[\frac{D_{q,\theta,\e}}{Dt}, \nabla\Delta^{-1}\mathbb P_{\approx \lambda_{q+1}}\right](T[W_{q+1}]\cdot\nabla \widetilde \theta_{q-1})\\
&+\nabla\Delta^{-1}\mathbb P_{\approx \lambda_{q+1}}\frac{D_{q,\theta,\e}}{Dt} (T[W_{q+1}]\cdot\nabla \widetilde \theta_{q-1}).
\end{split}
\end{equation}

The advective derivative of $R_{O,3}, ..., R_{O,6}$ can be estimated similarly and the details are omitted.

To estimate the material derivative of $R_M$, we first consider the term $(T[\widetilde \theta_{q-1}]- T[\widetilde \theta_{q-1,\e}])W_{q+1}$,
\begin{equation}\notag
\begin{split}
&\|(\partial_t+T[M_{q+1}]\cdot\nabla)\left( (T[\widetilde \theta_{q-1}]- T[\widetilde \theta_{q-1,\e}])W_{q+1}\right)\|_{C^0}\\
\lesssim&\ \|(\partial_t+T[M_{q,\e}]\cdot\nabla)\left( (T[\widetilde \theta_{q-1}]- T[\widetilde \theta_{q-1,\e}])\right)W_{q+1}\|_{C^0}\\
&+\|(\partial_t+T[M_{q,\e}]\cdot\nabla)W_{q+1}(T[\widetilde \theta_{q-1}]- T[\widetilde \theta_{q-1,\e}])\|_{C^0}\\
&+\|(T[M_{q}]-T[M_{q,\e}])\cdot\nabla \left( (T[\widetilde \theta_{q-1}]- T[\widetilde \theta_{q-1,\e}])\right)W_{q+1}\|_{C^0}\\
&+\|(T[M_{q}]-T[M_{q,\e}])\cdot\nabla W_{q+1} (T[\widetilde \theta_{q-1}]- T[\widetilde \theta_{q-1,\e}])\|_{C^0}\\
&+\|(T[W_{q+1}]\cdot\nabla) \left( (T[\widetilde \theta_{q-1}]- T[\widetilde \theta_{q-1,\e}])\right)W_{q+1}\|_{C^0}\\
&+\|(T[W_{q+1}]\cdot\nabla) W_{q+1}(T[\widetilde \theta_{q-1}]- T[\widetilde \theta_{q-1,\e}])\|_{C^0}.
\end{split}
\end{equation}
We see from Lemma \ref{le-est-adv} that the cost of $(\partial_t+T[M_{q,\e}]\cdot\nabla)$ is $\lambda_q\delta_{q-1}^{\frac12}$. Combining with Lemma \ref{le-molli} and estimate (\ref{est-correct-1}) we deduce 
\begin{equation}\notag
\begin{split}
&\|(\partial_t+T[M_{q,\e}]\cdot\nabla)\left( (T[\widetilde \theta_{q-1}]- T[\widetilde \theta_{q-1,\e}])\right)W_{q+1}\|_{C^0}\\
\lesssim&\ \|(\partial_t+T[M_{q,\e}]\cdot\nabla)\left( (T[\widetilde \theta_{q-1}]- T[\widetilde \theta_{q-1,\e}])\right)\|_{C^0}\|W_{q+1}\|_{C^0}\\
\lesssim&\ \lambda_q\delta_{q-1}^{\frac12} \lambda_{q}^{-1}\lambda_{q-1}\delta_{q-2}^{\frac12}\delta_{q}^{\frac12},
\end{split}
\end{equation}
\begin{equation}\notag
\begin{split}
&\|(\partial_t+T[M_{q,\e}]\cdot\nabla)W_{q+1} (T[\widetilde \theta_{q-1}]- T[\widetilde \theta_{q-1,\e}])\|_{C^0}\\
\lesssim&\ \|(\partial_t+T[M_{q,\e}]\cdot\nabla)W_{q+1}\|_{C^0}\|(T[\widetilde \theta_{q-1}]- T[\widetilde \theta_{q-1,\e}])\|_{C^0}\\
\lesssim&\ \lambda_q\delta_{q-1}^{\frac12}\delta_{q}^{\frac12} \lambda_{q}^{-1}\lambda_{q-1}\delta_{q-2}^{\frac12}.
\end{split}
\end{equation}
Using Lemma \ref{le-molli}, (\ref{induct-q-7}) and estimate (\ref{est-correct-1}),
\begin{equation}\notag
\begin{split}
&\|(T[M_{q}]-T[M_{q,\e}])\cdot\nabla \left( (T[\widetilde \theta_{q-1}]- T[\widetilde \theta_{q-1,\e}])\right)W_{q+1}\|_{C^0}\\
\lesssim&\ \|T[M_{q}]-T[M_{q,\e}]\|_{C^0}\|\nabla  (T[\widetilde \theta_{q-1}]- T[\widetilde \theta_{q-1,\e}])\|_{C^0}\|W_{q+1}\|_{C^0}\\
\lesssim&\  \lambda_{q+1}^{-1}\lambda_q\delta_{q-1}^{\frac12}\lambda_{q-1}\lambda_q^{-1}\lambda_{q-1}\delta_{q-2}^{\frac12}\delta_q^{\frac12},
\end{split}
\end{equation}
\begin{equation}\notag
\begin{split}
&\|(T[M_{q}]-T[M_{q,\e}])\cdot\nabla W_{q+1} (T[\widetilde \theta_{q-1}]- T[\widetilde \theta_{q-1,\e}])\|_{C^0}\\
\lesssim&\ \|T[M_{q}]-T[M_{q,\e}]\|_{C^0}\|T[\widetilde \theta_{q-1}]- T[\widetilde \theta_{q-1,\e}]\|_{C^0}\|\nabla W_{q+1}\|_{C^0}\\
\lesssim&\  \lambda_{q+1}^{-1}\lambda_q\delta_{q-1}^{\frac12}\lambda_{q}^{-1}\lambda_{q-1}\delta_{q-2}^{\frac12}\lambda_{q}\delta_q^{\frac12},
\end{split}
\end{equation}
Applying Lemma \ref{le-molli} and estimate (\ref{est-correct-1}),
\begin{equation}\notag
\begin{split}
&\|(T[W_{q+1}]\cdot\nabla) \left( (T[\widetilde \theta_{q-1}]- T[\widetilde \theta_{q-1,\e}])\right)W_{q+1}\|_{C^0}\\
\lesssim&\ \|T[W_{q+1}]\|_{C^0}\|\nabla(T[\widetilde \theta_{q-1}]- T[\widetilde \theta_{q-1,\e}])\|_{C^0}\|W_{q+1}\|_{C^0}\\
\lesssim&\ \delta_{q}^{\frac12}\lambda_{q-1} \lambda_{q}^{-1}\lambda_{q-1}\delta_{q-2}^{\frac12}\delta_{q}^{\frac12},
\end{split}
\end{equation}
\begin{equation}\notag
\begin{split}
&\|(T[W_{q+1}]\cdot\nabla) W_{q+1}(T[\widetilde \theta_{q-1}]- T[\widetilde \theta_{q-1,\e}])\|_{C^0}\\
\lesssim&\ \|T[W_{q+1}]\|_{C^0}\|T[\widetilde \theta_{q-1}]- T[\widetilde \theta_{q-1,\e}]\|_{C^0}\|\nabla W_{q+1}\|_{C^0}\\
\lesssim&\ \delta_{q}^{\frac12} \lambda_{q}^{-1}\lambda_{q-1}\delta_{q-2}^{\frac12}\lambda_{q}\delta_{q}^{\frac12}.
\end{split}
\end{equation}
Other terms in $R_M$ can be estimated analogously. 

\medskip

\textbf{H\"older estimates in time:} 
It follows from
\begin{equation}\notag
\partial_t W_{q+1}=D_{t,q} W_{q+1}-u_q\cdot\nabla W_{q+1}
\end{equation}
that
\[\|W_{q+1}\|_{C_t^1 C_x^0}\lesssim \delta_{q-1}^{\frac12}\lambda_{q+1}\delta_q^{\frac12},\]
and hence by interpolation
\[\|W_{q+1}\|_{C_t^\zeta C_x^0}\lesssim (\lambda_{q+1}\delta_{q-1}^{\frac12})^{\zeta}\delta_q^{\frac12}\sim \lambda_q^{(b-\frac{1}{2b}\beta)\zeta-\frac{1}{2}\beta}. \]
The $C^\zeta$ regularity in time is assured if
\[(b-\frac{1}{2b}\beta)\zeta-\frac{1}{2}\beta<0\]
i.e.
\[\zeta<\frac{\beta}{2b-\frac{1}{b}\beta}<\frac1{2d}\]
when choosing $b=1^+$ and $\beta<\frac{2}{2d+1}$.

\bigskip

\section{Proof of Theorem \ref{thm}}
\label{sec-proof}

Due to the presence of an external forcing in each equation of (\ref{pm-q}), there is redundancy to find an initial tuple $(P_0, M_0,  \bar R_0, \widetilde R_0)$ satisfying the system (\ref{pm-q}) and (\ref{Rq-decomp})-(\ref{induct-q-6}) at level 0, such that $M_0\not\equiv 0$. Applying the inductive Proposition \ref{prop-main} iteratively we obtain a sequence $\{(P_q, M_q, \bar R_q, \widetilde R_q\}$ satisfying  (\ref{pm-q}) and (\ref{Rq-decomp})-(\ref{induct-q-6}). Thanks to the estimates (\ref{induct-q-1}), (\ref{induct-q-3}) and (\ref{induct-q-5}), taking the limit $q\to \infty$ yields a limit solution $(P, M, \bar R, 0)$ of (\ref{pm1}) with $M\not\equiv 0$ and $f=\nabla\cdot \bar R$. Hence there are at least two solutions to (\ref{ase}) with external forcing $f=\nabla\cdot \bar R$. It follows from (\ref{Rt-q1}), (\ref{Rt-q2}) and the estimates in Proposition \ref{prop-main}  that $f\in C_t^0C_x^{2\alpha-1}$.



\end{document}